\input amstex
\input epsf

\hfuzz 35pt

\def\v{\vert}

\def\wwtilde#1{\widehat#1}
\def\Pr{\Bbb P}
\def\EE{\Bbb E}
\def\la{\langle}
\def\ra{\rangle}
\def\ol#1{\overline{#1}}
\def\O{\text{\rm O}}
\def\o{\text{\rm o}}
\def\g{\gamma}
\redefine\ell{l}
\redefine\Tilde{\widetilde}
\redefine\Re{\Cal R}
\def\v{\vert}
\def\vep{\varepsilon}
\documentstyle{grgppt01}
\redefine\qed{\hfill$\square$}
\pagewidth{125mm}
\pageheight{185mm}
\parindent=8mm
%\TagsOnRight
\magnification=1200
\def\section#1{\bigskip\goodbreak\line{{\bf\ignorespaces #1}\hfill}
  \nobreak\smallskip\nobreak\noindent\ignorespaces}

\def\subsection#1{\goodbreak\medskip
  \flushpar{\smc #1}\nobreak\smallskip \nobreak\flushpar\ignorespaces}

\def\capt#1#2{\baselineskip=10pt\botcaption{\baselineskip=10pt\eightpoint{\it Figure #1.}\quad
  #2}\endcaption\par}

\def\mletter#1#2#3{\hskip#2cm\lower#3cm\rlap{$#1$}\hskip-#2cm}
\def\lastletter#1#2#3{\hskip#2cm\lower#3cm\rlap{$#1$}\hskip-#2cm\vskip-#3cm}
\def\figure#1\par{\parindent=0pt
  \vbox{\baselineskip=0pt \lineskip=0pt
  \line{\hfil}
  #1}}

\catcode`\@=11
 \def\logo@{}
\catcode`\@=13
\def\today{\number\day
\space\ifcase\month\or
  January\or February\or March\or April\or May\or June\or
  July\or August\or September\or October\or November\or December\fi
  \space\number\year}
\topmatter
\title Random Electrical Networks on Complete Graphs II: Proofs
\endtitle
\author Geoffrey Grimmett and Harry Kesten 
\endauthor
\address
(GRG, current address) Statistical Laboratory, University of Cambridge, 
Wilberforce Road, Cambridge CB3 0WB, United Kingdom
\endaddress
\email g.r.grimmett{\@}statslab.cam.ac.uk
\endemail
\http http://www.statslab.cam.ac.uk/$\sim$grg/ \endhttp
\address
(HK) Mathematics Department, Cornell University,
Ithaca NY 14853, USA
\endaddress
\email kesten{\@}math.cornell.edu\endemail

\abstract
This paper contains the proofs of Theorems 2 and 3 of
the article entitled
{\it Random electrical networks on complete graphs\/},
written by the same authors and published in the
Journal of the London Mathematical Society, vol.\ 30 (1984), pp.\  171--192.
The current paper was written in 1983 but was
not published in a journal, although its existence
was announced in
the LMS paper. This \TeX\ version was created on
9 July 2001. It incorporates minor improvements to formatting and
punctuation, but no change has been made to the mathematics.

We study the effective electrical resistance of the complete
graph $K_{n+2}$ when each edge is allocated a random resistance.
These resistances are assumed independent with distribution
$\Pr(R=\infty)=1-n^{-1}\gamma(n)$, $\Pr(R\le x) = n^{-1}\gamma(n)F(x)$
for $0\le x < \infty$, where $F$ is a fixed distribution function
and $\gamma(n)\to\gamma\ge 0$ as $n\to\infty$.
The asymptotic effective resistance between two chosen vertices is 
identified in the two cases $\gamma\le 1$ and
$\gamma>1$, and the case $\gamma=\infty$ is considered. 
The analysis proceeds via detailed estimates based on
the theory of branching processes.
\endabstract
\endtopmatter
\document

\section{1. Introduction} In these notes we give complete proofs of
Theorems 2 and 3 and a further indication of the proof of Theorem 1 in
Grimmett and Kesten (1983). We use the same notation as in that paper and
we therefore repeat only the barest necessities.  $K_{n+2}$ denotes the
complete graph with $n+2$ vertices, which we label as $\{0,1,\dots
,n,\infty\}$. (See Bollob\'as (1979) for definition).  Each edge $e$ is
given a random resistance $R(e)$ with distribution
$$
\aligned
\Pr(R(e) \leq x) & = \frac{\g(n)}{n} F(x) \quad \text{for } 0 \leq x < \infty \\
\Pr(R(e)  = \infty) & = 1 - \frac{\g(n)}{n},
\endaligned\tag 1.1
$$
where $F$ is a fixed distribution function concentrated on $[0,\infty)$
and $\g(n)$ a sequence of numbers such that $0 \leq \g (n) \leq n$.  All
the resistances $R(e)$, $e \in K_{n+2}$, are assumed independent.  $R_n$
denotes the resulting (random) effective resistance in $K_{n+2}$ between
the vertices $0$ and $\infty$.  We shall prove the following result (the
numbering is taken from Grimmett and Kesten (1983)):

\proclaim{Theorem 2} If
$$
\lim_{n\to \infty} \g (n) = \g \leq 1 \tag 1.2
$$
then
$$
\lim_{n\to\infty} \Pr \big(R_n = \infty\big) = 1. \tag 1.3
$$
\endproclaim

To describe the limit distribution of $R_n$ when $\g (n) \to \g >1$ we
need a (one-type) Bienaym\'e--Galton--Watson process $\{Z_n\}_{n \geq 0}$
in which the offspring distribution is a Poisson distribution with mean
$\g$ and $Z_0=1$. (See Harris (1963) Ch. I; this book uses the more
traditional name Galton--Watson process for the branching process).  We
denote the random family tree of such a process by $T$ and label its root
by $\la 0\ra$, and the children in the $n$th generation of the individual
$\la i_1,\dots , i_{n-1}\ra$ (or $\la 0 \ra$ if $n=1$) in the $(n-1)$th
generation by $\la i_1,\dots , i_n\ra$ with $1 \leq i_n \leq N=N(i_1,\dots
, i_{n-1}):= \text{number of children of }\la i_1,\dots , i_{n-1}\ra$.  
Thus, not
all $\la i_1,\dots , i_n\ra$ with $i_1,\dots , i_n \geq 1$ occur as
vertices of $T$, but only those which correspond to individuals which are
actually born or ``realized''.  For more details, see Harris (1963), Ch.\ 
VI.2 or Jagers (1975), Ch.\ 1.2 or Grimmett and Kesten (1983).  
$T_0:=\{\la 0\ra\}$ is called the $0${\it th generation\/} of $T$, and
for $n \geq 1$ the collection of vertices $\la i_1,\dots , i_n\ra$ of $T$
is called the $n${\it th generation\/} of $T$, and denoted by $T_n$.
$\vert T_n\vert$, the cardinality of $T_n$, is just $Z_n$.  The subtree of
$T$ consisting of all vertices in $T_0 \cup T_1\cup \dots \cup T_n$
together with all the edges between these vertices is denoted by
$T_{[n]}$. We write $R(T_{[n]})$ for the resistance between $\la 0\ra$ and
$T_n$ in $T_{[n]}$ (formally this is defined by first identifying all
vertices in $T_n$ --- or shortcircuiting them --- and finding the resistance
between $\la 0 \ra$ and the single vertex obtained by this
identification).  Note that $R(T_{[n]})=\infty$ if and only if
$T_n=\emptyset$, i.e., if and only if the branching process is extinct by
the $n$th generation.  The limit
$$
R(T): = \lim_{n\to \infty} R(T_{[n]})
$$
always exists by the monotonicity property (2.6) below.

\proclaim{Theorem 3}
If
$$
\lim \g(n) = \g>1, \tag 1.4
$$
then, as $n\to\infty$, the distribution of $R_n$ converges to that of
$R'(\g)+R''(\g)$, where $R'(\g)$ and $R''(\g)$ are independent random
variables, each with the distribution of $R(T)$, defined above, with the
mean of the Poisson offspring distribution equal to the value $\g$ given
by (1.4).  In particular the atom at $\infty$ of the limit distribution of
$R_n$ equals $2q(\g)-q^2 (\g)$ $(<1)$, where $q(\g)$ is the extinction
probability of the branching process with Poisson offspring distribution
with mean $\g$ ($q(\g)$ is the smaller solution of the equation $q=\exp (-\g
(1-q))$). 
\endproclaim

\remark{Remark}
Theorems 2 and 3 show that the value one is a critical value for $\g$.  
If $\g \leq 1$, then all the mass in the distribution of $R_n$ escapes to
$\infty$ as $n\to \infty$, while for $\g > 1$, $R_n$ has a limit
distribution which puts some mass on $[0,\infty)$ (but also has an atom at
$\infty$).  This is closely related to the threshold phenomenon for the
spread of epidemics discussed by von Bahr and Martin-L\"of (1980,
especially Section 5). Both the von Bahr and Martin-L\"of paper and ours
rest on the existence of imbedded random trees which behave like the
family trees of branching processes (see next paragraph).

The idea of the proofs was already explained in Grimmett and Kesten
(1983).  It consists in looking at the graphs of vertices of $K_{n+2}$
which are connected to $0$ and $\infty$, respectively, by paths of finite
resistance (such paths will be called {\it conducting paths\/} in the
sequel).  It will be shown that these graphs resemble two independent
trees, each with the same distribution as $T$, given above.  In addition,
it will be shown that there are with high probability a great many edges
with finite resistance joining pairs of vertices, one from each tree,
which are far away from $0$ and $\infty$, respectively.  There are enough
of these interconnections to make $R_n$ nearly equal to the sum of the
resistances of these two trees (one connected to $0$ and one to $\infty$).  
All this will be done first under the assumption that $R(e) \geq \vep> 0$
for all $e$ with probability one.  The next section is largely devoted to
proving the continuity property in Proposition 1 which allows us to let
$\vep$ go to $0$ afterwards.  This continuity property needs proof because
$T$ may be infinite.  In finite networks continuity of the effective
resistance between two vertices, as a function of the resistances of the
individual edges, is comparatively easy (see Kesten (1982) Ch.\  11).
\endremark

\section{2. Preliminaries} 
Some standard ways to combine resistances were
already discussed in Section 2 of Grimmett and Kesten (1983).  We need
some (known) extensions of these rules, especially for the case where
individual edges may have zero resistance.

Let $G$ be a finite connected graph and $A_0$ and $A_1$ two disjoint sets
of vertices of $G$.  Assume that each edge $e$ has been assigned a
resistance, to be denoted by $R(e)$.  To find the resistance between $A_0$
and $A_1$ in the network of edges of $G$ one first identifies all vertices
in $A_0$ ($A_1$) as a single vertex, $\widehat{A}_0$ ($\widehat{A}_1$) say.  
This is equivalent to setting $R(e)$ equal to zero, whenever both
endpoints of $e$ lie in the same $A_i$.  We shall also identify as one
vertex any maximal class of vertices of $G$ which is already
shortcircuited, i.e., any maximal class $\widehat{A}=\{v_1,\dots , v_m\}$
such that for any $v_i, v_j$ in $\widehat{A}$ there exist
$v_{i_{1}},\dots, v_{i_{r}}$ and edges $e_\ell$ between $v_{i_{\ell}}$ and
$v_{i_{\ell+1}}$, $\ell=0,\dots,r$, with $v_{i_{0}}=v_i$,
$v_{i_{r+1}}=v_j$, and $R(e_\ell)=0$, $\ell=0,\dots, r$.  Let $\widehat{G}$
be the network resulting from these identifications.  The resistance
between $A_0$ and $A_1$ in $G$ is defined as the resistance between
$\widehat{A}_0$ and $\widehat{A}_1$ in $\widehat{G}$.  To compute this
resistance one introduces the potential function $V(\widehat{v})$ (with
$\widehat{v}$ running through the vertices of $\widehat{G}$) with the
boundary values $0$ on $\widehat{A}_0$ and $1$ on $\widehat{A}_1$ (to
produce these boundary values physically, one has to connect
$\widehat{A}_0$ and $\widehat{A}_1$ to a voltage source external to the
network).  $V(\cdot)$ is determined by Kirchhoff's laws:
$$
V(\widehat{v}) = \bigg\{\sum \frac{1}{R(e)}\bigg\}^{-1} 
\sum\frac{V(\widehat{w}(e))}{R(e)}, \quad \widehat{v}\neq \widehat{A}_0, \ 
\widehat{A}_1.
\tag 2.1
$$
The sums in (2.1) run over all edges $e$ of $G$ with one endpoint in the
class corresponding to $\widehat{v}$ and the other endpoint outside this
class; the class of the endpoint of $e$ outside $\widehat{v}$ is denoted
by $\widehat{w}(e)$.  Note that any $R(e)$ appearing in (2.1) is strictly
positive by our choice of the classes $\widehat{v}$ and $\widehat{w}$ (see
Kesten (1982) Ch.\ 11 for more details).

We shall frequently appeal to the following probabilistic interpretation
of $V(\cdot)$ (see Doyle and Snell (1982), Griffeath and Liggett (1982)).  
Consider a Markov chain $\{X_\nu\}$ on $\widehat{G}$ with transition
probability
$$
P(\widehat{v},\widehat{w}) = \bigg\{\sum_v\frac{1}{R(e)}\bigg\}^{-1}
\sum_{v,w} \frac{1}{R(e)}, \quad \widehat{v} \neq \widehat{w},\tag 2.2
$$
where $\sum_v$ runs over all edges $e$ of $G$ with one endpoint inside and
one endpoint outside $\widehat{v}$, while $\sum_{v,w}$ runs only over
those edges with one endpoint in $\widehat{v}$ and the other in
$\widehat{w}$.  Then
$$
V(\widehat{v}) = \Pr\{X_\cdot\text{ visits } \widehat{A}_1 
\text{ before } \widehat{A}_0\mid X_0 = \widehat{v}\}.\tag 2.3
$$

The resistance between $\widehat{A}_0$ and $\widehat{A}_1$ is given by
$$
\bigg\{\sum \frac{V(\widehat{w}(e))}{R(e)} \bigg\}^{-1},\tag 2.4
$$
with the sum in (2.4) running over all edges $e$ of $G$ with one endpoint
in $A_0$, and the other endpoint in any class disjoint from $A_0$.  (In
(2.4) the class of this other endpoint is denoted by $\widehat{w}(e)$;
$\widehat{w}(e)$ varies with $e$.) See Kesten (1982) Ch.\ 11.  The
probability interpretation (2.3) together with (2.4) gives a probabilistic
meaning to resistance as well.  The reader should note that the above
simplifies if $R(e)>0$ for all $e$ and if $A_0$ and $A_1$ consist of
single vertices. In this case $\widehat{G}=G$.

Very intuitive is the following monotonicity property.  Let $G$, $A_0$ and
$A_1$ be as above, and let $\{R'(e)\}$, $\{R''(e)\}$ be two assignments of
resistances to the edges of $G$.  Denote the corresponding resistances
between $A_0$ and $A_1$ in $G$ by $R'(A_0,A_1)$ and $R''(A_0, A_1)$,
respectively.  Then
$$
R'(e) \leq R''(e) \quad \text{for all } e \tag 2.5
$$
implies
$$
R'(A_0, A_1) \leq  R''(A_0, A_1).\tag 2.6
$$

Unfortunately the proof is not all that simple (see Griffeath and Liggett
(1982), Doyle and Snell (1982), and for the case when $R'(e)$ and $R''(e)$
may take the values $0$, $\infty$ see Kesten (1982)).  Note that this
monotonicity property states in particular that shortcircuiting some
vertices, or insertion of additional edges (no matter what their
resistance is) can only decrease the resistance between $A_0$ and $A_1$;
also removal of any edges can only increase the latter resistance.

For the remainder of this section $Z_0=1, Z_1, Z_2,\dots$ is any
Bienaym\'e--Galton--Watson branching process with the mean number $\g$ of
offspring per individual strictly greater than 1, but finite.  That is to
say,
$$
1 < \g: = \EE Z_1 <\infty. \tag 2.7
$$
It is not assumed that the offspring distribution is a Poisson
distribution.  $q$ denotes the extinction probability:
$$
q = \Pr\{Z_n = 0 \text{ eventually}\}.\tag 2.8
$$
It is well known (see Harris (1963) Theorem I.6.1) that under (2.7)
$$
q<1. \tag 2.9
$$
We write $f(z)$ for the generating function of the offspring distributions:
$$
f(z) = \sum^\infty_{n=0} \Pr\{Z_1 = n\} z^n, \quad \v z\v \leq 1.\tag 2.10
$$
The convexity of $f$ on $[0,1]$ and (2.9) imply that $f'(q)<1$ (see Harris 
(1963) Fig.\ 1 and proof of Theorem I.8.4).  We can therefore find an 
$\vep_0>0$ such that
$$
0<q + 2 \vep_0 <1, \ f'(q+2\vep_0) + 2\vep_0 < 1. \tag 2.11
$$

$T$ will be the family tree of the branching process as in Section 1.  
Also $T_n$ and $T_{[n]}$ are as in Section 1 and each edge of $T$ is 
given a resistance such that the $\{R(e) : e \in T\}$  are independent, 
and all have the same 
distribution $F$. 
 As in Section 1 we define
$$
R(T) = \lim_{n\to\infty} R(T_{[n]}) = \lim_{n\to\infty}
\{\text{resistance between 0 and } T_n \text{ in } T_{[n]}\}.
$$
Since we can think of $R(T_{[n]})$ as the resistance between 0 and 
$T_{n+1}$ in $T_{[n+1]}$ when all vertices in $T_n \cup T_{n+1}$ are 
shortcircuited, it follows from the monotonicity property (2.6) 
that $R(T_{[n]}) \leq R(T_{[n+1]})$ so that the limit $R(T)$ is well 
defined so long as we allow it to take the value $\infty$. We set
$$
R^\vep(e) = R(e) + \vep,
$$
and in general, for any resistance $R(*)$ calculated as a function of the
$R(e)$, we denote by $R^\vep(*)$ the corresponding resistance when $R(e)$
is everywhere replaced by $R^\vep (e)$.  In particular $R^\vep (T)$ is the
resistance of the family tree when each edge resistance is increased by
$\vep$.

\proclaim{Proposition 1}
$$
\lim_{\vep \downarrow 0} R^\vep (T) = R(T) \quad \text{w.p.1}.\tag 2.12
$$
\endproclaim

The proof will be broken down into several lemmas. If $\la i_1,\dots, i_n\ra $
is a vertex of $T$ then we write $T(i_1,\dots, i_n)$ for the subtree of
$T$ whose vertices are $\la i_1,\dots, i_n\ra $ and all its descendants.  I.e.,
the vertex set of $T(i_1,\dots, i_n)$ is
$$
\la i_1,\dots, i_n\ra \bigcup \Bigl\{\la i_1,\dots, i_n, \ j_1,\dots, 
j_\ell \ra  \in T: \ell \geq 1, 
j_1,\dots, j_\ell \in \{1,2,\dots \}\Bigr\}.
$$
Two vertices of $T(i_1,\dots, i_n)$ have an edge of $T(i_1,\dots, i_n)$ 
between them if and only if they are connected by an edge in $T$.  
From the branching property, it follows that, conditionally 
on $\la i_1,\dots, i_n\ra \in T$, the distribution
of $T(i_1,\dots, i_n)$ (as a random graph) is the same 
as the original distribution of $T$.  More generally, given 
that $\la i_1,\dots, i_m\ra \in T$ and 
$1 \leq n_1 < n_2 <\dots<n_\ell \leq m-1$, 
$j_r \neq i_{n_{r+1}}$ and $\la i_1,\dots, i_{n_{r}}$, $j_r\ra \in T$,
$$
\aligned
&T(i_1,\dots, i_{n_{r}}, j_r), \ r = 1,\dots, \ell, \text{ are conditionally
independent, each}\\
&\text{with the same distribution as $T$,
and the edges of these}\\
&\text{trees have independent
resistances, all with distribution } F.
\endaligned
\tag 2.13
$$
We also introduce the notation $T^j(i_1,\dots, i_n)$ for the subtree of
$T$ whose vertices are $\la i_1,\dots, i_n\ra $, $\la i_1,\dots, i_n, j\ra $ and all
descendants of $\la i_1,\dots, i_n, j\ra $. Note that $\la i_1,\dots, i_n\ra $ only has
the one descendant $\la i_1,\dots, i_n, j\ra $ in this tree.  Similarly
$T^j=T^j(0)$ has as vertices $\la 0\ra , \la j\ra $ and the descendants of $\la j\ra $.  We
only use this notation if $\la i_1,\dots, i_n,j\ra \in T$ (respectively $\la j\ra  \in
T$ if $n=0$). Given that $\la i_1,\dots, i_n,j\ra \in T$ the distribution of
$T^j(i_1,\dots , i_n)$ is the same as the conditional distribution of $T$,
given $\v T_1\v =1$. \footnote{$\v A\v$ denotes the number of vertices
in $A$.}) There is also an independence property for several $T^{j_{r}}
(i_1,\dots, i_{n_{r}})$ analogous to (2.13), given that $\la i_1,\dots,
i_{n_{r}}, j_r\ra \in T$.

\proclaim{Lemma 1}
Let $\Cal P$ be a property of
rooted labeled trees with resistances assigned to their
edges.  If $\vep_0$ satisfies (2.11), and if
$$ 
\Pr\{T^j \text{\rm\ does not have property } \Cal P \mid
 \la j\ra  \in T_1\} \leq q
+ 2 \vep_0, \tag 2.14
$$
then
$$
\align
\Pr\{&\text{\rm for each infinite path } \la i_1, i_2,\ldots \ra  \text{\rm\ in } T
\text{\rm\ there exist infinitely}\tag "\rlap{(2.15)}"\\
&\text{\rm many } n \text{\rm\ and integers }
j_{n+1} \neq i_{n+1} \text{\rm\ such that } \la i_1,\dots, i_n, j_{n+1} \ra  \in T\\
&\text{\rm and such that } T^{j_{n+1}} (i_1,\dots, i_n) \text{\rm\ has property
} \Cal P \} = 1. 
\endalign
$$
Moreover, if $\Gamma (i_1, i_2,\dots , i_n)$ denotes the number of $k$, $2
\leq k \leq n$ for which there exists a $j_k \neq i_k$ such that
$T^{j_{k}} (i_1,\dots, i_{k-1})$ has property $\Cal P$, then there exist
constants $0 < C_1, C_2 < \infty$ such that for $n \geq 2$
$$
\Pr\Bigl\{T_n \neq \emptyset 
\text{\rm\ and } \min_{\la i_1,\dots, i_n\ra \in T_n} \Gamma
(i_1,\dots, i_n) \leq C_1 n\Bigr\} \leq e^{-C_{2}n} .\tag 2.16
$$
\endproclaim

\demo{Proof}  Assume (2.14).  By virtue of (2.13) and its analogue for several 
$T^{j_{r}} (i_1,\dots, i_{n_{r}})$ it suffices for (2.15) to prove that
$$\align
\lim_{n\to\infty} \sum\Sb i_1\ge 1,\dots,\\i_n\ge1\endSb
\Pr \Bigl\{
&\la i_1,\dots , i_n\ra  \in T \text{ but there does not 
exist a } k \leq n\tag 2.17\\
&\text{and } j_k \neq i_k,
\text{ such that } \la i_1, \dots , i_{k-1}, j_k\ra  \in T \\
&\text{and such that } T^{j_{k}} (i_1,\dots , i_{k-1}) \text{ has 
 property } \Cal P\Bigr\} = 0.
\endalign
$$
To see this, note that (2.17) says that if $\la i_1, i_2,\ldots \ra $ is an
infinite path in $T$, then there is with probability one at least one $n$
and $j_{n+1} \neq i_{n+1}$ for which $T^{j_{n+1}} (i_1,\dots , i_n)$ has
property $\Cal P$. But $T^{j_{n+1}} (i_1,\dots, i_n)$ and $T^{i_{n+1}}
(i_1,\dots , i_n)$ are independent, and (2.17) again applies to
$T^{i_{n+1}} (i_1,\dots, i_n)$, so that with probability one there exists
a further $m>n$ and $j_{m+1}\ne i_{m+1}$ such that $T^{i_{m+1}} (i_1,\dots,
i_m)$ also has property $\Cal P$ etc.

To prove (2.17) note that $\la i_1,\dots, i_n \ra \in T$ if and only if $\la 0\ra $
has $\ell_1 \geq i_1$ children, $\dots$, $\la i_1,\dots, i_r\ra $ has
$\ell_{r+1} \geq i_{r+1}$ children, $r=0,\dots, n-1$, for some integers
$\ell_{r+1}$.  The probability of this event for given $\ell_{r+1} \geq
i_{r+1}$ is
$$
p_{\ell_{1}} \cdot p_{\ell_{2}}\cdots p_{\ell_{n}} ,
$$
where $p_\ell = \Pr\{Z_1 = \ell\}$.  Given that the above event occurs with
prescribed $\ell_1,\dots, \ell_n$, the trees $T^{j_{r+1}} (i_1,\dots,
i_r)$, $r=0,\dots, n-1$, $j_{r+1} = 1,\dots, \ell_{r+1}$, $j_{r+1}\neq
i_{r+1}$ with their resistances are all independent, and each has the
conditional distribution of $T^j$ and its resistances, given $\la j\ra \in T_1$.  
Therefore the conditional probability that none of these trees has
property $\Cal P$ is at most
$$
(q+ 2\vep_0)^{\sum^{n}_{r=1} (\ell_r-1)} 
\tag 2.18
$$
(by virtue of (2.14)).  The summand of (2.17) is therefore at most
$$
\sum_{\ell_1 \geq i_1, \dots, \ell_n \geq i_n}
 p_{\ell_{1}}\cdots p_{\ell_{n}} (q+2\vep_0)^{\sum^{n}_{r=1} (\ell_{r}-1)} .
$$
Finally, the sum in (2.17) is at most
$$
\align
\sum_{\ell_1\geq 1,\dots, \ell_n \geq 1} \sum\Sb 1 \leq i_1\le \ell_1\\
\vdots\\ 1\le i_n\le \ell_n\endSb
 p_{\ell_{1}} \cdots p_{\ell_{n}} (q+2 \vep_0)^{\sum^n_{r=1} 
(\ell_r-1)} 
& = \biggl\{ \sum_{\ell \geq 1} \ell p_\ell (q+2\vep_0)^{\ell-1}\biggr\}^n  
\tag 2.19\\
&= \{f' (q+2\vep_0)\}^n .
\endalign
$$
By (2.11) the last member of (2.19) tends to zero as $n\to\infty$, so that
(2.15) holds.

A slight strengthening of the above argument leads to (2.16).  For every
$\theta \geq 0$ and $\ell_{r+1} \geq i_{r+1}$,
$$
\aligned
\Pr\{&\la i_1,\dots, i_n\ra  \in T_n 
  \text{ and } \la i_1,\dots , i_r\ra  \text{ has }
\ell_{r+1} \text{ children in } T,\\
&\hskip3cm  r=0,\dots, n-1, \text{ but } \Gamma 
(i_1,\dots, i_n) \leq C_1 n \}\\
\leq{}&p_{\ell_{1}} \cdots p_{\ell_{n}} e^{\theta C_{1}n} \EE \{e^{-\theta 
\Gamma (i_{1},\dots, i_{n})} \mid \la i_1,\dots, i_r\ra  \text{ has } 
\ell_{r+1} \text{ children in } T, \ r=0,\dots, n-1\}.
\endaligned
\tag 2.20
$$
Denote the left hand side of (2.14) temporarily by $\rho$.  Analogously to
(2.18) we then obtain
$$\align
\EE\{&e^{-\theta \Gamma (i_{1},\dots, i_{n})}\mid\la i_1,\dots, i_r\ra  
\text{ has } \ell_{r+1} \text{ children, }
 0 \leq r \leq n-1\}\tag "\rlap{(2.21)}"\\
 &= \prod^{n-1}_{r=1} \{\rho^{\ell_{r+1}-1} + 
(1-\rho^{\ell_{r+1}-1}) e^{-\theta}\}\\
& \leq \prod^{n-1}_{r=1} \{(q+2\vep_0)^{\ell_{r+1}-1} + 
(1-(q+2\vep_0^{\ell_{r+1}-1}) e^{-\theta}\}.
\endalign
$$
Moreover
$$
\align
\sum^\infty_{i=1} \sum_{\ell \geq i} {}&p_\ell \{(q+2\vep_0)^{\ell-1} + 
(1-(q+2\vep_0)^{\ell-1}) e^{-\theta}\} \tag 2.22\\
= & \sum^\infty_{\ell=1} \ell p_\ell \{(q+2\epsilon_0)^{\ell-1} 
(1-e^{-\theta}) + e^{-\theta}\}\\
= & (1-e^{-\theta}) f' (q+2\vep_0) + e^{-\theta} \gamma .
\endalign
$$
Now choose $\theta>0$ so large that
$$
(1-e^{-\theta}) f' (q+2\vep_0) + e^{-\theta} \gamma \leq f' (q+2\vep_0) 
+ \vep_0 \leq 1 - \vep_0
$$
(see (2.11)).  Since the left hand side of (2.16) is bounded by the 
sum over $\ell_1 \geq 1, \dots, \ell_n \geq 1$, $i_1 \leq \ell_1,\dots, 
i_n \leq \ell_n$ of the left hand side of (2.20), we obtain from 
(2.20)--(2.22)
$$
\align
&\Pr\{\min_{\la i_{1},\dots, i_{n} \ra\in T} \Gamma (i_1,\dots, i_n) 
\leq C_1 n \}\\
&\hskip.5cm\leq e^{\theta C_{1}n} \sum_{i_{1} \geq 1,\ldots, i_{n} \geq 1} 
\quad\sum_{\ell_{1} \geq i_{1},\dots, \ell_{n} \geq i_{n}} p_{\ell_{1}} 
\cdots p_{\ell_{1}}\\ 
&\hskip4cm \times\prod^{n-1}_{r=1} \{(q+2\vep_0)^{\ell_{r+1}-1} +   
(1-(q+2\vep_0)^{\ell_{r+1}-1}) e^{-\theta}\}\\
&\hskip.5cm \leq e^{\theta C_{1} n} (1-\vep_0)^{n-1} .
\endalign
$$
(2.16) follows if we take $C_1$ small enough.
\qed\enddemo

\proclaim{Lemma 2}
For all $\vep >0$
$$
\Pr\{R^\vep (T) = \infty\} = \Pr\{R(T) = \infty\} = q .
$$
\endproclaim

\demo{Proof}
We prove the second equality.  This proof works for any choice of $F$, and
thus implies the first equality as well, since adding $\vep$ to each
resistance has the same effect on the distribution of $R(T)$ as changing
$F(x)$ to $F(x-\vep)$.

By the rules for combining resistances (see Grimmett and Kesten (1983)
Section 2; also Figure 1 below),
$$\align
\{R(T)\}^{-1} & = \sum_{1 \leq j \leq Z_{1}} \{R(T^{j})\}^{-1}\tag 2.23\\
& = \sum_{1 \leq j \leq Z_{1}} \{R(e(j)) + R(T(j))\}^{-1}, 
\endalign
$$
where $e(j)$ denotes the edge between $\la 0\ra $ and $\la j\ra $. $(T(j)$ is defined
just after Proposition 1).

\topinsert
\figure   
\mletter{T^1}{1.7}{3.9}
\mletter{T^2}{4.1}{1}
\mletter{T^3}{10.3}{4}
\mletter{e(1)}{4.7}{4.1}
\mletter{e(2)}{6.1}{3.1}
\lastletter{e(3)}{7}{3.9}
\centerline{\epsfxsize=8cm\epsffile{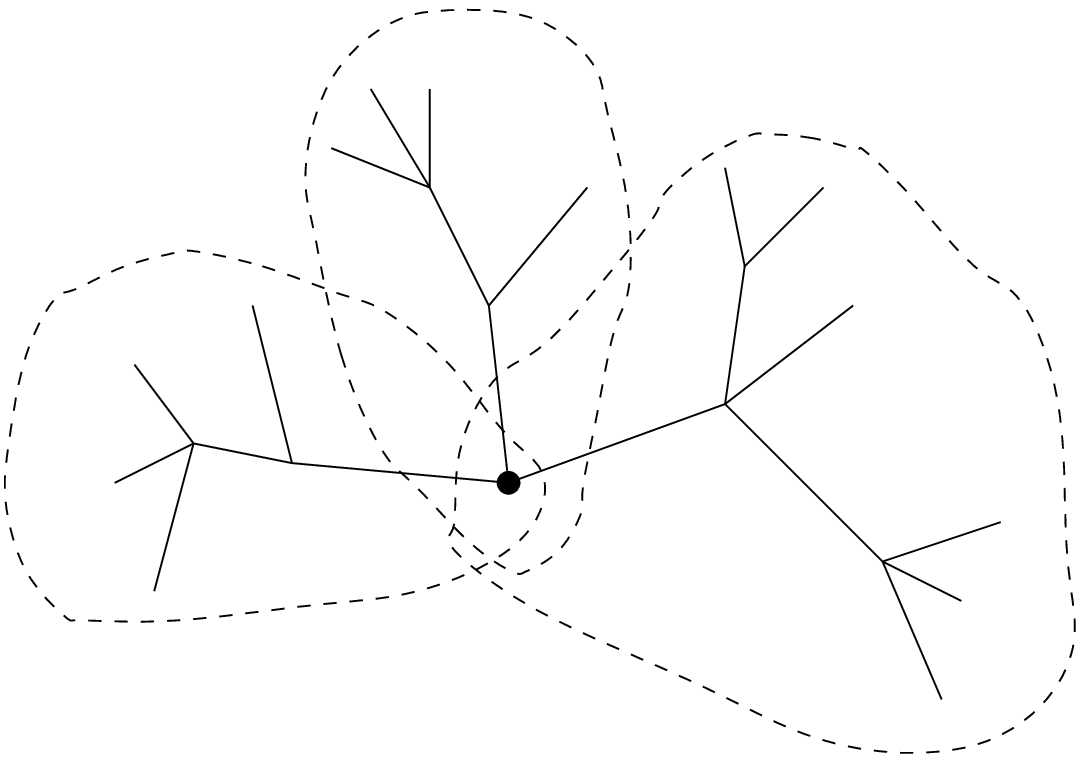}}

\capt{1}{A picture of $T$ when $Z_1=3$; the root of $T$, $\la 0\ra $, is
represented by the heavy dot.  The subtrees $T^j$ form parallel
resistances between $\la 0\ra $ and $\infty$.  
In $T^j$, $e(j)$ and $T(j)$ are in
series.}
\endinsert

In particular $R(T)=\infty$ if and only if $Z_1=\v T_1\v = 0$ or
$R(T(j))=\infty$ for each $j \leq Z_1$. In view of (2.13) this implies
$$\align
\Pr\{R(T)=\infty\} & = \sum_{n \geq 0} \Pr\{Z_1 = n\} \Pr\{R(T(j)) = 
\infty \text{ for }
1 \leq j \leq n\} \\
& = \sum_{n \geq 0} \Pr\{Z_1 = n\} (\Pr \{R(T)=\infty\})^n.
\endalign
$$
Thus $\Pr\{R(T)=\infty\}$ is a solution of the equation
$$
x=f(x).\tag 2.24
$$
As is well known (see Harris (1963) proof of Theorem I.6.1) the only
solutions of (2.24) in $[0,1]$ are $q$ and $1$.  It therefore suffices to
show that
$$
\Pr \{R(T) = \infty\}<1.\tag 2.25
$$
To this end we first choose a constant $K$ such that
$$
\gamma F(K) >1 . \tag 2.26
$$
We next consider the subtree $\widetilde{T}$ of $T$ whose vertices are the
vertices of $T$ connected to $\la 0\ra $ by a path all of whose edges have a
resistance not exceeding $K$.  Of course two vertices of $\widetilde{T}$
are connected by an edge of $\widetilde{T}$ if and only if they are
connected by an edge of $T$.  We write $\widetilde{T}_n$ (respectively
$\widetilde{T}_{[n]}$) for the part of $\widetilde{T}$ which belongs to
$T_n$ (respectively $T_{[n]}$).  The children in $\widetilde{T}$ of a
vertex $\la x\ra $ in $\widetilde{T}$ are 
precisely those connected to $\la x\ra $ by
an edge of resistance not exceeding $K$. $\widetilde{T}$ is therefore the
family tree of a branching process $\{\widetilde{Z}_n\}$ with offspring
distributions
$$
\widetilde{p}_m:= \Pr\{\widetilde{Z}_1=m\} = \sum_{n \ge m} \Pr\{Z_1 = n\}
\binom nm  F^m (K)(1-F(K))^{n-m},
$$
and mean number of offspring
$$
\sum_{m \geq 0} m \widetilde{p}_m = 
\sum_{n \geq 0} n \Pr\{Z_1 = n\} F(K) = \gamma F(K) >1.
\tag 2.27
$$
(See (2.26)).  Thus the $\widetilde{Z}$ process is supercritical also, and
$$
\widetilde{q}:= \Pr \{\widetilde{Z}_n = 0  \text{ eventually}\} < 1.
$$
From Theorem I.6.2 of Harris (1963) we conclude that for each $m \geq 1$
$$
\Pr \{\widetilde{Z}_k \geq m\} \to 1 - \widetilde{q} >0 
\quad\text{ as } k \to \infty .
$$
We can therefore fix $\kappa$ such that
$$
\sum_{n \geq 2} \Pr\{\widetilde{Z}_\kappa = n\} 
\{1- \big(\tfrac{1}{2} + \tfrac{\widetilde{q}}{2}\big)^n - n 
\big(\tfrac{1}{2} - \tfrac{\widetilde{q}}{2}\big)\big(\tfrac{1}{2} + 
\tfrac{\widetilde{q}}{2}^{n-1}\} \geq \tfrac{3}{4} 
(1-\widetilde{q}).\tag 2.28
$$
With $\kappa$ fixed in this way we define
$$\align
g(x) & = \sum_{n \geq 2} \Pr (\widetilde{Z}_\kappa = n) \{1-(1-x)^n - nx(1-x)^{n-1}\} \\
     & = \sum_{n\geq 2} \Pr(\widetilde{Z}_\kappa = n) \sum^n_{j=2} 
\binom nj x^j (1-x)^{n-j}, \quad 0 \leq x \leq 1, \\
L & = 2 \kappa K 
\endalign
$$
and
$$
\alpha_n = \Pr(R(\widetilde{T}_{[\kappa n]}) \leq L) \quad (\alpha_0 = 1),
$$
where $R(\widetilde{T}_{[n]})$ is the resistance between 
$\la 0\ra $ and $\widetilde{T}_n$ in $\widetilde{T}_{[n]}$. We claim that
$$
\alpha_n \geq g (\alpha_{n-1}),\quad n \geq 1 . \tag 2.29
$$

Before proving (2.29) we show that it quickly implies the lemma.  Clearly
$g(x)$ is non-decreasing and continuous on $[0,1]$ and (2.28) states that
$g(\frac12(1-\widetilde{q})) \geq \tfrac{3}{4}
(1-\widetilde{q}) > \tfrac{1}{2} (1-\widetilde{q})$, while $g(x) \leq 1$.  
Therefore
$$
r: = \max \{x\in [0,1]: g(x) \geq x\} >\tfrac{1}{2} (1-\widetilde{q}) > 0.
$$
Thus $r=g(r) \leq g(x) < x$ for $x \in (r,1]$, which together with (2.29) 
and $\alpha_0 =1$ implies (see Figure 2)
$$
\alpha_n \geq g(\alpha_{n-1}) \geq g(g(\dots (g(1))\dots )) \to r 
\quad\text{as } n\to \infty .
$$
Thus
$$
\lim_{n\to\infty} \Pr (R(\widetilde{T}_{[\kappa n]}) \leq L) \geq r.
$$
By the monotonicity property (2.6) $R(T_{[n]}) \leq
R(\widetilde{T}_{[n]})$ so that also $\Pr (R(T)\leq L)\geq r$.  Thus
(2.29) will imply (2.25) and the lemma.

\topinsert
\figure
\mletter{y}{3.9}{.5}
\mletter{y=x}{8.6}{.9}
\mletter{y=g(x)}{5}{2.6}
\mletter{0}{3.9}{5.3}
\mletter{r}{6.6}{5.3}
\mletter{1}{8.1}{5.3}
\lastletter{x}{9}{5.5}
\centerline{\epsfxsize=6cm\epsffile{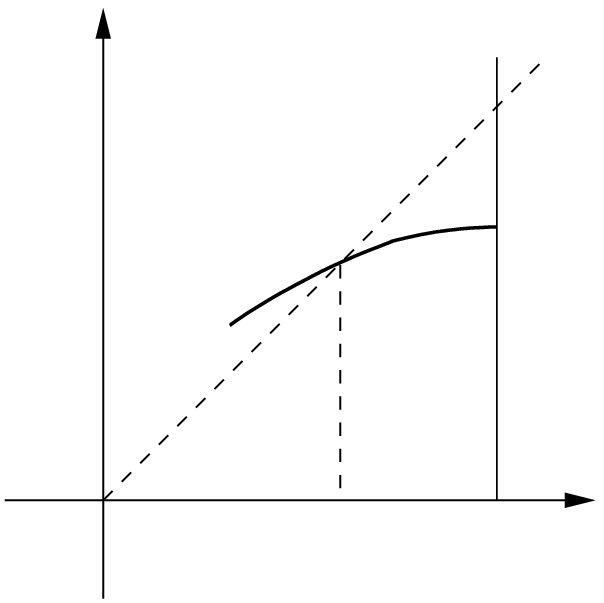}}

\capt{2}{}
\endinsert

Now we prove (2.29).  Consider the event that $\widetilde{T}_\kappa$
contains at least two distinct individuals $\la i_1,\dots, i_\kappa\ra $ and
$\la j_1,\dots, j_\kappa\ra $ such that the two resistances, between 
$\la i_1,\dots,
i_\kappa\ra $ and $\widetilde{T}_{\kappa n}$, and between $\la j_1,\dots,
j_\kappa\ra $ and $\widetilde{T}_{\kappa n}$, in the tree of the $\kappa
(n-1)$ generations of descendants of $\la i_1,\dots, i_\kappa\ra $ are both at most
$L$.  By virtue of (2.13) the probability of this event is
$$
\align
&\sum_{n\geq 2} \Pr(\widetilde{Z}_\kappa = n) \sum^n_{j=2} 
\binom nj (\Pr \{R(\widetilde{T}_{[\kappa (n-1)]} 
\leq L)\})^j 
(1-P \{R(\widetilde{T}_{[\kappa (n-1)]} \leq L)\})^{n-j}\tag 2.30\\
&\hskip9cm = g(\alpha_{n-1}).
\endalign
$$
Moreover, if $i_1 = j_1,\dots, i_s = j_s$, $i_{s+1} \neq j_{s+1}$, 
then the network consisting of the edges from $\la 0\ra $ to $\la i_1,\dots, i_s\ra $ 
and the two parallel connections from $\la i_1,\dots, i_s\ra $ to 
$\widetilde{T}_{\kappa n}$ via $\la i_1,\dots, i_\kappa\ra $ and 
the tree of its descendants and via $\la j_1,\dots, j_\kappa\ra  = 
\la i_1,\dots, i_s$, $j_{s+1},\dots, j_\kappa\ra $ 
and the tree of its descendants is at most (see Figure 3)
$$
sK + \tfrac{1}{2} \{(\kappa -s) K + L \} \leq L;
$$
recall that each edge of $\widetilde{T}$ has resistance $\leq
K$, and $L=2 \kappa K$. Thus $R(T_{[\kappa n]}) \leq L$ whenever
$\la i_1,\dots, i_\kappa\ra $ and $\la j_1,\dots, j_\kappa\ra $ exist as above.  
Consequently (2.29) follows from the value in (2.30) for the probability
of the existence of such $\la i_1,\dots, i_\kappa\ra $ and
$\la j_1,\dots, j_\kappa\ra $.
\qed\enddemo

\topinsert
\figure
\mletter{\widetilde T_{\kappa n}}{9}{1.3}
\mletter{\la i_1,\ldots,i_\kappa\ra}{3.3}{4.3}
\mletter{\la j_1,\ldots,j_\kappa\ra}{8}{4.5}
\mletter{\la i_1,\ldots,i_s\ra}{6.7}{5.9}
\lastletter{\la 0\ra}{6.7}{7.5}
\centerline{\epsfxsize=7cm\epsffile{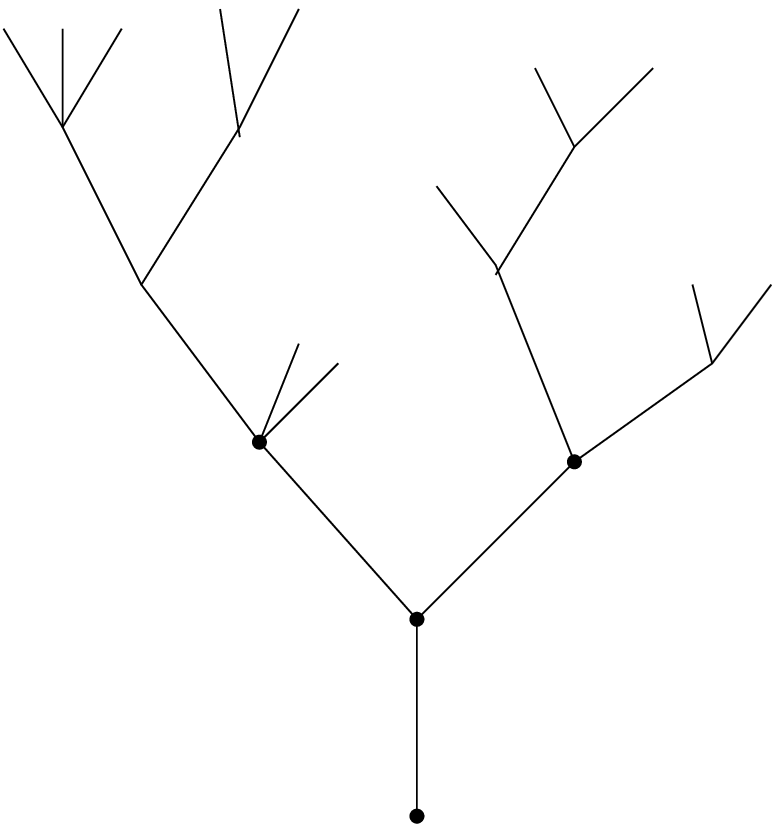}}

\capt{3}{}
\endinsert

We need some more notation for the next lemma.  Note that this lemma does
not involve random quantities.  Let $t$ be a rooted labeled tree, with
root $\la 0\ra $ and vertices labeled $\la i_1,\dots, i_n\ra $, just as described for
$T$ in the introduction.  Assume that to each edge $e$ of $t$, a resistance
$r(e)<\infty$ has been assigned.  If $\la x\ra =\la i_1,\dots, i_n\ra $ is a vertex of
$t$, then set
$$\align
\rho(x)&= \{\text{resistance of the unique path in } t
\text{ from } \la 0\ra  \text{ to } \la x\ra\} \tag "\rlap{(2.31)}"\\
&= \sum^n_{k=1} \{\text{resistance of edge between } \la i_1,\dots, i_{k-1}\ra  
\text{ and } \la i_1,\dots, i_k\ra \}
\endalign 
$$
(for $k=1$, $\la i_1,\dots, i_{k-1}\ra =\la 0\ra $). In accordance with previous
notation we write $r(t(x))$ for the resistance of the tree $t(x)$
consisting of $\la x\ra $ and all its descendants and connecting edges between
them.  As before, $r(t(x))$ is the limit (as $m\to\infty$) of the
resistance in $t(x)$, between $\la x\ra $ and its descendants in the $(n+m)$th
generation of $t$, i.e., between $\la x\ra =\la i_1,\dots,i_n\ra $ and the collection
of vertices of the form $\la i_1,\dots, i_n$, $j_1,\dots, j_m\ra $.  Finally, if
$\la y\ra  \in t$ its {\it equivalence class\/} 
$\la \widehat{y}\ra $ is the collection
of vertices of $t$ connected to $\la y\ra $ by paths of zero resistance (see the
beginning of this section).

\proclaim{Lemma 3}
Let $t$ be a rooted labeled tree as above.  Let $\{X_\nu\}_{\nu \geq 0}$
be a Markov chain with state space the collection of equivalence classes
$\{\la \widehat{y}\ra :\la y\ra  \in t\}$ and transition probability matrix
$$
P(\la y\ra , \la z\ra ) = \biggl\{\sum_y \frac{1}{r(e)}\biggr\}^{-1} 
\sum\nolimits_{y,z} \frac{1}{r(e)}, \quad \la \widehat{y}\ra  \neq \la \widehat{z}\ra  ,
\tag 2.32
$$
where $\sum_y$ runs over all edges $e$ of $t$ with one endpoint in
$\la \widehat{y}\ra $ and one outside $\la \widehat{y}\ra $, while $\sum_{y,z}$ runs
only over those edges with one endpoint in $\la \widehat{y}\ra $ and the other
in $\la \widehat{z}\ra $ (compare (2.2)). Then\footnote{The second member of
(2.33) stands of course for $\lim_{n\to\infty} \Pr (X.$ reaches the $n$th
generation of $t$ before 
$\la \widehat{0}\ra \mid X_0=\la \widehat{x}\ra \}.$ We shall
use the more intuitive expression of (2.33) even in proofs without
formally going through taking the limit as $n\to\infty$.}
$$\align
\Pr(X. \text{\rm\ never reaches } \la \widehat{0}\ra 
   \mid X_0 = \la \widehat{x}\ra )
&= \Pr (X. \text{\rm\ reaches } \infty \text{\rm\ before } 
\la \widehat{0}\ra \mid X_0 = \la \widehat{x}\ra )\tag 2.33
\\
&\geq  \frac{\rho(x)}{\rho(x) + r (t(x))} 
, \quad \la \widehat{x}\ra  \neq \la \widehat{0}\ra ,
\endalign
$$
provided $\rho(x) + r(t(x))>0$.
\endproclaim

\demo{Proof} For simplicity we only consider the case where all $r(e)>0$
so that each equivalence class $\la \widehat{y}\ra $ consists of one vertex
only.  Since the identification of vertices in the same equivalence class
turns $t$ into another tree \footnote{Actually this is an abuse of
terminology since we may have multiple edges between a pair
$\la \widehat{y}\ra $ and $\la \widehat{z}\ra $.} with vertices $\la \widehat{y}\ra $, it
is not hard to extend the argument to the general case where $r(e)=0$ is
allowed.

\topinsert
\figure
\mletter{t(x)}{7.7}{.5}
\mletter{\la i_1,\ldots,i_n\ra}{6.2}{2.46}
\mletter{\la i_1,i_2\ra}{6.2}{4.66}
\mletter{\la i_1\ra}{6.2}{5.4}
\lastletter{\la 0\ra}{6.2}{6.1}
\centerline{\epsfysize=6cm\epsffile{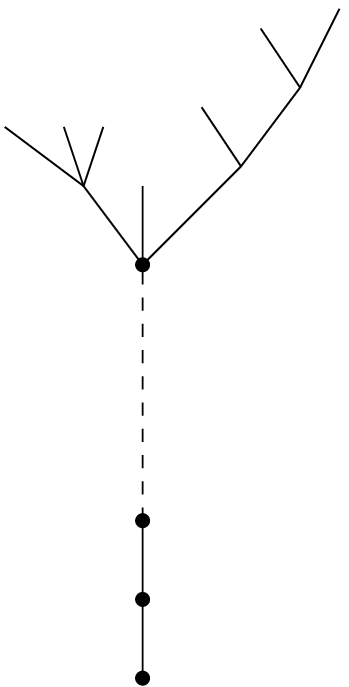}}

\capt{4}{The tree $t_*$.}
\endinsert

Now let $\la x\ra  = \la i_1,\dots, i_n\ra $ and consider the tree $t_*$, whose
vertices are only $\la 0\ra , \la i_1\ra , 
\la i_1, i_2\ra ,\ldots, \la i_1,\dots, i_n\ra =\la x\ra $ and
all the descendants of $\la x\ra $ (see Figure 4).  In this tree the
resistance
between $\la 0\ra $ and $\la x\ra $ is precisely $\rho(x)$, and the resistance between
$\la x\ra $ and $\infty$ is precisely $r(t(x))$. Let $\{X_{*_{\nu}}\}_{\nu\geq
0}$ be the Markov chain on $t_*$ which is analogous to $X_\nu$, and has
transition probability matrix 
$$ 
P_* (\la y\ra , \la z\ra ) = \biggl\{\sum_{y}{}\kern-2pt\lower7pt\hbox{$*$}
\frac{1}{r(e)}\biggr\}^{-1} \biggl\{\sum_{yz}{}\kern-2pt\lower7pt\hbox{$*$} \frac{1}{r(e)}\biggr\}, 
\quad \la z\ra  \neq \la y\ra , 
$$
where $\sum_{*y}$ (respectively $\sum_{*yz}$) runs over all edges of
$t_*$ with one endpoint at $\la y\ra $ and the other endpoint different from
$\la y\ra $ (respectively at $\la z\ra $).  From the relation between potentials and
hitting probabilities (see (2.3)) we see that $\Pr_*\{X_*.$ reaches
$\infty$ before $\la 0\ra \mid 
X_{*0} = \la x\ra \rbrace =V_*(\la x\ra )$, where $V_*$ is the
potential in the network $t_*$ when $\la 0\ra $ (respectively $\infty$) is given
potential zero (respectively one).  But the resistance between $\la x\ra $ and
$\la 0\ra $ in $t_*$ equals $\rho(x)$, and the resistance between $\la x\ra $ and
$\infty$ equals $r(t(x))$.  Standard computations based on Kirchhoff's
laws now show $V_* (\la x\ra )=(\rho(x)+r(t(x))^{-1}\rho (x)$, or equivalently
$$
\Pr_*\{X_{*.} \text{ reaches } \infty \text{ before } \la 0\ra  \mid X_{*0} =
\la x\ra \}
= \frac{\rho (x)}{\rho(x) + r(t(x))}.
\tag 2.34 
$$
The right hand sides of (2.33) and (2.34) are the same.  On the other
hand, $t_*$ is formed by removing from $t$ all descendants of
$\la i_1,\dots, i_k\ra $ other than $\la i_1,\dots, i_{k+1}\ra $ for $k=0,1,\dots,
n-1$. Now consider the successive times $\nu_i$ at which $X_.$ visits
$t_*$ at a different point than at the last visit to $t_*$.  Formally,
$\nu_0=0$, 
$$ 
\nu_{i+1} = \inf \{\nu > \nu_i: X_\nu \in t_*, X_\nu \neq
X_{\nu_{i}}\rbrace. 
$$ 
If $X_{\nu_{j}} = \la i_1,\dots, i_k\ra \in t_*$ with
$k< n$, then $X_.$ may make several excursions into $\bigcup_{j\neq i_{k+1}}
t(i_1,\dots, i_k, j)$ before it visits a point $\la x\ra \in t_*$ which differs
from $\la i_1,\dots, i_k\ra $.  However,
 $X_.$ can reach such a vertex in $t_*$
in one step only from $\la i_1,\dots, i_k\ra $. It must first return from an
excursion into
 $t(i_1,\dots, i_k, j)$, $j \neq i_{k+1}$, to $\la i_1,\dots,i_k\ra $ before it 
can reach $\la i_1,\dots, i_{k-1}\ra $ or $\la i_1,\dots, i_{k+1}\ra $. 
From this it is 
not hard to see that, given $X_0=\la x\ra \in t_*$ and $\nu_\ell <\infty$, the 
distribution\footnote{Correction added in 2001:
This sentence should not assert equidistribution, but instead that
$\Pr(X_{\nu_i}=\la x_i\ra,\ 0\le i\le l\mid X_0=\la x\ra)
\le \Pr_*(X_{*i}=\la x_i\ra,\ 0\le i\le l\mid X_{*0}=\la x\ra)$.}
of $X_0, X_{\nu_{1}},\dots , X_{\nu_{\ell}}$ is the same as that
of $X_{*0}, X_{*1},\dots, X_{*\ell}$. \footnote {$\{X_{*\nu}\}$ is
almost the imbedded chain of $\{X_\nu\}$ on $t_{*}$, but only observed at
times when the state changes.}
From this it is practically obvious that
$$
\align
&\Pr \{X_. \text{ reaches } \la 0\ra  \text{ at some time}  
  \mid X_0 = \la x\ra \}\\
&\hskip3cm= \Pr \{X_{\nu_{i}} = \la 0\ra  \text{ for some } i  \mid X_0 = \la x\ra \}\\
&\hskip3cm\leq \Pr_* \{X_{*\nu} = \la 0\ra  \text{ at some time}  \mid X_{*0} = \la x\ra \}.
\endalign
$$
Consequently
$$
\Pr \{X_.\text{ never reaches } \la 0\ra   \mid  X_0 = \la x\ra \}
\geq \Pr_* \{X_{*.} \text{ never reaches } \la 0\ra   \mid  X_{*0} = \la x\ra  \}.
$$
In view of (2.34) this is equivalent to (2.33).
\qed\enddemo

We now prove Proposition 1 in two lemmas, separating the cases 
$\gamma F(0) \leq 1$ and $\gamma F(0)> 1$,

\proclaim{Lemma 4}
If 
$$
\gamma F(0) \leq 1, 
$$
then (2.12) holds.
\endproclaim

\demo{Proof}
First define $\widetilde{Z}_n$ and $\widetilde{T}$ as in Lemma 2 with $K=0$.
Thus
$\widetilde{Z}_n$ is the branching process of nodes which are connected to 
$\la 0\ra $ by paths of zero resistance.  As in (2.27), the mean number of 
offspring per individual for this branching process is $\gamma F(0)$.  
But now $\gamma F(0) \leq 1$ so that the $\widetilde{Z}_n$ process dies out
eventually w.p.1. (see Harris (1963) Theorem I.6.1).  This means
that a.s.\ there is no infinite path in $T$ ---
the family tree of $Z_n$ --- all
of whose edges have zero resistance.

For $\la x\ra =\la i_1,\dots, i_n\ra  \in T$ define $\rho(x)$ as in (2.31) (with
$t$ replaced by $T$). Recall that the superscript in $R^\epsilon(\cdot)$
indicates that the resistances of all edges have been increased by $\epsilon$.
Now set
$$
R^+ (T(i_1,\dots, i_n)) = \lim_{\epsilon \downarrow 0} R^\epsilon(T(i_1,
\dots, i_n)).
$$
This limit exists by the monotonicity property (2.6).  The principal
estimate which we need is as follows: with probability 1, for
each infinite sequence $\{i_k\}_{k \geq 1}$ such that $\la i_1,\dots, i_n\ra 
\in T$ for all $n$ one has
$$
\liminf_{n\to\infty} \frac{R^+(T(i_1,\dots, i_n))}{\rho(\la i_1,\dots, i_n\ra )} 
= 0.\tag 2.35
$$
Write $e(i_1,\dots, i_n)$ for the edge between $\la i_1,\dots, i_{n-1}\ra $ and
$\la i_1,\dots, i_n\ra $ if $n>1$, and $e(i_1)$ for the edge between $\la 0\ra $
and $\la i_1\ra $. (2.35) is fairly easy for paths in $T$ which satisfy
$$
\sum^\infty_{k=1} R(e(i_1,\dots, i_k))=\infty . \tag 2.36
$$
For such a path we first choose $\epsilon_0>0$ such that (2.11) holds. By
Lemma 2 we can then find an $L<\infty$ such that
$$
\Pr \{R^1(T) >L\} \leq q + \epsilon_0.
$$
Then
$$\align
\Pr \{R^1(T^j)> 2L  \mid  \la j\ra  \in T_1\}
&\leq \Pr \{R^1 (e(j))> L\} + \Pr \{R^1 (T(j))> L \mid \la j\ra \in T_1\}\\
&\leq \Pr \{R^1 (e(1))>L\} + q + \epsilon_0,
\endalign
$$
since $R^1(T^j)=R^1(e(j))+R^1(T(j))$ whenever $\la j\ra \in T_1$
(compare (2.23)). Without loss of generality we can therefore choose $L$
large enough such that
$$
\Pr\lbrace R^1(T^j)>2L \mid \la j\ra \in T_1\rbrace\le q+2\varepsilon_0.\tag 2.37
$$
We now apply Lemma 1 with the following choice for $\Cal P$: $T$ has
property $\Cal P$ if and only if $R^1(T)\le 2L$. Then by (2.15), w.p.1, for
each infinite path $\la i_1,i_2,\dots \ra $ in $T$ there exist infinitely many
$n$ and $j_{n+1}\ne i_{n+1}$ with
$$
R^1(T^{j_{n+1}}(i_1,\dots,i_n))\le 2L.\tag 2.38
$$ 
but $T^{j_{n+1}}(i_1,\dots,i_n)$ is a subtree of $T(i_1,\dots,i_n)$ 
so that (2.38) and the monotonicity property (2.6) imply
$$
R^1(T(i_1,\dots,i_n))\le R^1(T^{j_{n+1}}(i_1,\dots,i_n))\le 2L.
$$
Moreover, by virtue of (2.36)
$$
\rho (\la i_1,\dots ,i_n\ra )=\sum\limits_{k=1}^n R(e(i_1,\dots,i_k))\to\infty,
$$
so that (2.35) holds under the condition (2.36) 
(note that $R^+\le R^1$ by (2.6) again).

\topinsert
\figure
\mletter{T^{j_1}(i_1,\dots,i_{n_{1}})}{1.5}{0}
\mletter{ T^{j_2}(i_1,\dots,i_{n_{2}})}{5.3}{0}
\mletter{T^{j_k}(i_1,\dots,i_{n_{k}})}{9.4}{0}
\mletter{\la i_1,\dots,i_N\ra}{.6}{4.2}
\mletter{\la i_1,\dots,i_{n_1}\ra}{2.6}{4.2}
\mletter{\la i_1,\dots,i_{n_2}\ra}{5.6}{4.2}
\lastletter{\la i_1,\dots,i_{n_k}\ra}{9.5}{4.2}
\line{\hfil}
\centerline{\epsfxsize=10cm\epsffile{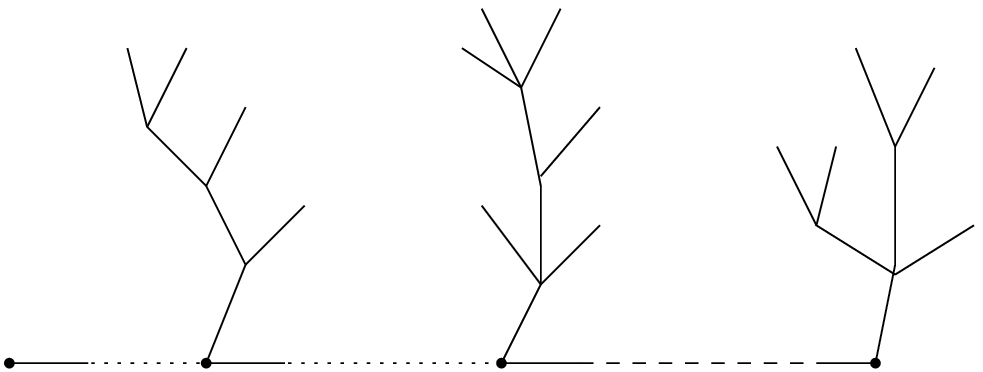}}

\capt{5}{}
\endinsert

If (2.36) fails for some path $i_1,i_2,\dots,$ then 
$$
\sum\limits_{k=1}^\infty R(e(i_1,\dots,i_k))<\infty\tag 2.39
$$
for this path. Such paths can actually have positive probability; 
see Bramson (1978). In this case we can find a.s.\ for each $\delta>0$ an
 $N$ such that 
$$
\sum\limits_{k=N+1}^\infty R(e(i_1,\dots,i_k))
\le\delta\rho (\la i_1,\dots,i_N\ra )
=\delta \sum\limits_{k=1}^N R(e(i_1,\dots,i_k))
\tag 2.40 
$$
(since a.s.\ the last sum is strictly positive for large $N$ by the 
first paragraph of the proof). Also, by the preceding argument we 
can find a.s.\ $N<n_1<n_2<\cdots$ and $j_\ell\ne i_{n_\ell +1}$ such that
$$
R^1(T^{j_\ell}(i_1,\dots,i_{n_\ell}))\le 2L.\tag 2.41
$$
In this case $T(i_1,\dots,i_N)$ contains the tree consisting of the path 
$\la i_1,\dots,i_N\ra$, $\la i_1,\dots,i_{N+1}\ra$,
$\dots$, $\la i_1,\dots,i_{n_k}\ra $
together with the trees $T^{j_\ell}(i_1,\dots,i_{n_\ell})$, $\ell
=1,\dots,k$ which are attached to this path. A schematic diagram of this
graph is given in Figure 5.
Since each $T^{j_\ell}(i_1,\dots,i_{n_{\ell}})$ satisfies (2.41) we also have 
$$
R^+(T^{j_\ell}(i_1,\dots,i_{n_{\ell}}))\le 2L,\quad j=1,\dots,k.
$$
Therefore $R^+(T(i_1,\dots,i_N))$ is at most equal to the resistance
between $v_0$ and $\infty$ in the network of Figure 6, where the
resistance of the edge between $v_i$ and $\infty$ is $2L$, and
the edge between $v_{\ell -1}$ and $v_\ell$ has resistance
$$
\beta_\ell :=\sum\limits_{r=n_{\ell-1}+1}^{n_\ell} R(e(i_1,\dots,i_r)) 
\quad (\text{with }n_0=N).
$$
A simple inductive argument (add a vertex on the left and a resistance
$2L$ between $v_0$ and $\infty )$ shows that the resistance between $v_0$
and $\infty$ is at most
$$
\align 
\beta_1 +\dots +\beta_k+\frac{2L}{k}&\le 
\sum\limits_{r=N+1}^\infty R(e(i_1,\dots,i_r))+\frac{2L}{k}
\tag "\rlap{(2.42)}"\\
&\le\delta\rho (\la i_1,\dots,i_N\ra )+\frac{2L}{k}\quad\text{ (see (2.40))}.
\endalign
$$
Since this estimate can be proved for each $k$, it follows that for each
$\delta >0$ we can find an $N$ with
$$
R^+(T(i_1,\dots,i_N))\le 2\delta\rho(\la i_1,\dots,i_N\ra ).
$$
This proves (2.35) in all cases.

\topinsert
\figure
\mletter{\infty}{5}{0}
\mletter{2L}{3.2}{3}
\mletter{2L}{4.2}{3}
\mletter{2L}{5.2}{3}
\mletter{2L}{7.4}{3}
\mletter{2L}{9.6}{3}
\mletter{v_0}{2.1}{4.1}
\mletter{\beta_1}{2.7}{4.1}
\mletter{v_1}{3.3}{4.1}
\mletter{\beta_2}{4}{4.1}
\mletter{v_2}{4.6}{4.1}
\mletter{v_3}{5.9}{4.1}
\mletter{v_{k-1}}{8.8}{4.1}
\mletter{\beta_{k-1}}{9.5}{4.1}
\lastletter{v_k}{10.3}{4.1}
\line{\hfil}
\centerline{\epsfxsize=8cm\epsffile{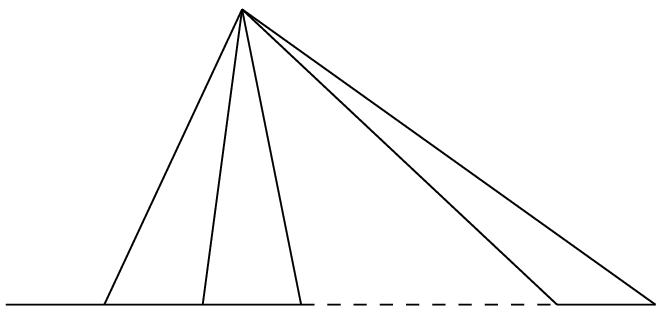}}

\capt{6}{The number next to an edge gives 
the resistance of that edge.}
\endinsert

We now prove (2.12) from (2.35). Fix $\delta >0$. We claim that there
exists a finite collection ${\Cal C}={\Cal C}(\delta )$ of vertices of
$T$ with the following properties:
$$
\alignat 2
&\text{Each $\la i_1,\dots,i_n\ra \in {\Cal C}$ satisfies: }
&&T(i_1,\dots,i_n)\text{ is infinite},\tag 2.43\\
&&&\rho (\la i_1,\dots,i_n\ra )>0\text{ and}\\
&&&R^+(T(i_1,\dots,i_n))\le\delta\rho (\la i_1,\dots,i_n\ra ).\\
&&&\llap{\text{Any path from $\la 0\ra$  to $\infty$ in }}
T  \text{ contains
some vertex in $\Cal C$.}\tag 2.44\\
\endalignat
$$
Note that the first paragraph of the proof and (2.35) show that w.p.1 
each infinite path $\la 0,i_1,i_2,\ldots \ra $ in $T$ contains a first 
vertex $\la i_1,\dots,i_n\ra $ with the properties listed in (2.43). 
Take for ${\Cal C}$ the collection of all these first vertices 
$\la i_1,\ldots,i_n\ra $ obtainable in this way as $\la 0,i_1,i_2,\dots \ra $ 
varies over the infinite paths in $T$. This ${\Cal C}$ has the properties 
(2.43), (2.44) and we merely have to verify that ${\Cal C}$ is finite. 
But, if ${\Cal C}$ were infinite, then there would have to exist an 
infinite sequence of infinite paths $\la 0,i_1^k,i_2^k,\ldots\ra \in T$, 
$k=1,2,\dots$, such that, for all $j\le k$, one of the properties listed 
in (2.43) fails for $\la i_1^k,\dots,i_j^k\ra $. By a diagonal 
selection we could then find an infinite path 
$\la 0,i_1,i_2,\dots\ra \in T$ such that, for each 
$j$, $\la i_1,\dots,i_j\ra $ lacks one of the properties in (2.43). 
Since we already saw that no such sequence exists ${\Cal C}$ must be finite.

Since ${\Cal C}$ is finite we can find an $\varepsilon_1>0$ such that 
for each $\varepsilon\le\varepsilon_1$ and each 
$\la i_1,\dots,i_n\ra \in {\Cal C}$
$$
R^\varepsilon (T(i_1,\dots,i_n))\le 2\delta\rho (\la i_1,\dots,i_n\ra )
\le 2\delta\rho^\varepsilon (\la i_1,\dots,i_n\ra ).
\tag 2.45  
$$ 
 Fix any $\varepsilon\le\varepsilon_1$ and let $X_\nu$ be the Markov chain
defined in Lemma 3 for $t=T$ and $r(e)=R^\varepsilon (e)$. Thus all
equivalence classes consist of a single vertex now. Then (see (2.4))
$$
\lbrace R^\varepsilon (T)\rbrace^{-1}=\sum\limits_{i\in T_1}
\frac{V(\la i\ra )}{R^\varepsilon (e(i))},
$$
where $V(\la i\ra )$ is the potential of $\la i\ra $ when $\la 0\ra $ is given the potential
zero and \footnote{As with $R^\varepsilon (T)$ this has to be
interpreted as a statement about a limit. Give $0$ potential 0 and $T_N$
potential $1$ and let $N\to\infty$. The potential $V(\la i\ra )$ is decreasing
in $N$ and hence has a limit, as one can see for instance from
(2.3).} $\infty$ is given potential one (by means of an external voltage
source). Note that $\la 0\ra \notin {\Cal C}$ by construction. It
follows therefore from (2.44) that if $X_0=\la i\ra \in T_1$, then
$X_\nu$ cannot reach $\infty$ without passing through ${\Cal C}$.
Consequently,
$$ 
\align 
V(\la i\ra )&=\Pr\lbrace X.\text{ reaches }\infty\text{ before
}\la 0\ra  \mid X_0=\la i\ra \rbrace\\ 
&=\sum_{\la x\ra \in {\Cal C}}\Pr\lbrace
X.\text{ reaches }{\Cal C}\text{ before }\la 0\ra \text{ and hits
}{\Cal C}\text{ first in }\la x\ra  \mid X_0=\la i\ra \rbrace\\
&\hskip3cm\times \Pr\lbrace X_.\text{
reaches } 
\infty\text{ before }\la 0\ra  \mid X_0=\la x\ra \rbrace. 
\endalign
$$
By Lemma 3 and (2.45) we therefore have
$$ 
\align 
&\lbrace R^\varepsilon
(T)\rbrace^{-1}\ge\frac{1}{1+2\delta}\sum_{i\in
T_1}\frac{1}{R^\varepsilon (e(i))}\tag "\rlap{(2.46)}"\\
&\times\sum_{\la x\ra \in {\Cal
C}}\Pr\lbrace X.\text{ reaches }{\Cal C}\text{ before }\la 0\ra \text{ and
 hits }{\Cal C}\text{ first in }\la x\ra  \mid X_0=\la i\ra \rbrace.
\endalign
$$

Let $\ol{T}$ be the finite subtree of $T$ which contains $\la 0\ra $ and all
vertices $\la i_1,\dots,i_k\ra $ which have no predecessor $\la i_1,\dots,i_\ell\ra $
with $\ell<k$ belonging to ${\Cal C}$. Then it follows again by (2.3) and
(2.4) that the double sum in the right hand side of (2.46) is just the
reciprocal of the resistance between $\la 0\ra $ and ${\Cal C}$ in $\ol{T}$,
still when each edge $e$ is given the resistance $R^\varepsilon (e)$.
However, $\ol{T}$ is a finite network, so that resistances in $\ol{T}$
are continuous in $\varepsilon$ as $\varepsilon$ tends to zero. This
finally gives
$$
\align 
R^+(T)=\lim_{\varepsilon\downarrow 0}R^\varepsilon (T)&\le
 (1+2\delta )\times\lbrace\text{resistance between }\la 0\ra 
\text{ and }{\Cal C}\text{ in }\ol{T}\rbrace\\
&\le (1+2\delta )R(T).
\endalign
$$
The last inequality is again an easy consequence of the monotonicity
property (2.6). Since, also by (2.6), $R^+(T)\ge R(T)$ we now obtain
(2.12) by letting $\delta\downarrow 0$.
\phantom{xxxxxxxxx}\qed\enddemo

\proclaim {Lemma 5}
If
$$
\gamma F(0)>1
$$
then (2.12) holds.
\endproclaim

\demo {Proof}
The proof of (2.35) remains valid for any path $i_1,i_2,\dots$ with

$$
\sum\limits_{k=1}^\infty R(e(i_1,\dots,i_k))=\lim_{n\to\infty} 
\rho (\la i_1,\dots,i_n\ra )>0.
\tag 2.47
$$
Consequently, the proof of Lemma 4 goes through unchanged on any
realization $T$ of the family tree and its resistances for which (2.47)
holds for all infinite paths in the tree, or equivalently for any
realization in which the equivalence class $\la \widehat 0\ra $ is finite. The only
new complication is that for $\gamma F(0)>1$ the branching process
$\widetilde{Z}_n$ of Lemma 4 is supercritical and has a positive probability
of never dying out (Harris (1963), Theorem I.6.1). If $\widetilde{Z}_n$ does
not die out, then $\la \widehat 0\ra $ is infinite. In this case $\la 0\ra $ is connected
to $\infty$ by a path of zero resistance and $R(T)=0$. We merely have to
show that also $R^+(T)=0$ a.s.\ on the set of realizations which contain a
path of zero resistance between $\la 0\ra $ and $\infty$. This, however, follows
also by the arguments of Lemma 4. With probability one for any path
$i_1,i_2,\dots$ of zero resistance there exist infinitely many $n$ and
$j_{n+1}\ne i_{n+1}$ for which (2.38) holds. The middle expression in
(2.42) with $N=0$ is therefore still an upper bound for $R^+(T)$ whenever
such a path $i_1,i_2,\dots$ of zero resistance exists. But this says
$R^+(T)\le 2L/k$ for any $k$ or $R^+(T)=0$ a.s.\ on the set of trees for
which (2.47) fails for some path. Thus the lemma holds in all cases.
\qed\enddemo 

Lemmas 4 and 5 together prove Proposition 1.

We shall need an improved version of Lemma 3 for the case where the
resistance assigned to $e$ is at least $\varepsilon >0$. (In particular
this will hold if we assign the value $R^\varepsilon (e)$ to $e$.) In this
situation the equivalence classes considered in Lemma 3 all consist of one
vertex only. For each realization $t$ and $r(\cdot )$ of $T$ and its
resistances we can therefore consider the Markov chain $\lbrace
X_\nu^m\rbrace_{\nu\ge 0}=\lbrace X_\nu^m(t_{[m]},r(\cdot ))\rbrace$ with
state space the vertices in $t_{[m]}$ and transition probability matrix
$$
P(\la y\ra ,\la z\ra )=
\cases \biggl\{\displaystyle\sum_y\dfrac{1}{r(e)}\biggr\}^{-1}\dfrac{1}{r(y,z)}
&\text{if }\la y\ra \text{ and }\la z\ra \text{ are adjacent in }t_{[m]},\\
0&\text{ otherwise}.\endcases
\tag 2.48
$$
Here $r(y,z)$ is the resistance of the edge between $y$ and $z$, and $\sum_y$ is the sum over all edges $e$ in $t_{[m]}$ which are incident to $y$. For $s<m$ let 
$$
\align
{\Cal A}_s=\lbrace \la y\ra \in t_s:{}&\la y\ra \text{ has at least one descendant in}\\
&t_m,\text{ the $m$th generation of }t_{[m]}\rbrace.
\endalign
$$  
Further, for $\la x\ra \in t_{s-1}$ set
$$
\pi (x,t_{[m]},r,s)=\Pr\lbrace X_.^m\text{ reaches }\la 0\ra \text{ before }
\text{it reaches }{\Cal A}_s \mid X_0^m=\la x\ra \rbrace
$$
and
$$
\Pi (t_{[m]},r)=\Pi (t_{[m]},r,s)=\max_{\la x\ra \in 
t_{s-1}}\pi(x,t_{[m]},r,s).\tag 2.49
$$
Thus, $\Pi$ measures the probability for $X_\nu^m$ to go from $t_{s-1}$
through $t_{[m]}$ to the root of $t_{[m]}$ without passing through ${\Cal
A}_s$. When $t_{[m]}$ and $r$ are taken random again, then $\Pi$ is also a
random variable.

\proclaim {Lemma 6}
If $R(e)\ge\varepsilon$ w.p.1, then there exist constants
$0<C_1,C_2,L<\infty$ (independent of $\varepsilon$, $s$ and $m$) such that
$$
\Pr\left\lbrace{\Cal A}_s\ne\emptyset \text{\rm\ and }\Pi (T_{[m]},R,s)>
\left(\frac{2L}{2L+\varepsilon}\right)^{C_{1}s}\right\rbrace
\le e^{-C_{2}s}\quad\text{for }s\ge 2.\tag 2.50
$$
\endproclaim

\demo {Proof}
Let $\la x\ra =\la i_1,\dots,i_{s-1}\ra \in t_{s-1}$. 
Then, if $X_0^m=\la x\ra$,
$X_\nu^m$ cannot reach $\la 0\ra $ without passing 
through $\la i_1,\dots,i_k\ra $ for
each $1\le k\le s-1$. Let $\tau_k$ be the first time $X_.^m$ reaches
$\la i_1,\dots,i_k\ra$ ($\tau_{s-1}=0$). 
Then for $\tau_k\le\nu<\tau_{k-1}$,
$X_\nu$ takes only values in the collection of descendants of
$\la i_1,\dots,i_{k-1}\ra $ in $T_{[m]}$ (this is the collection of vertices
$\la i_1,\dots,i_{k-1},j_k,\dots,j_\ell\ra$, $k\le\ell\le m)$. 
These are vertices
of $T(\la i_1,\dots,i_{k-1}\ra )$, and therefore between $\tau_k$ and
$\tau_{k-1}$, $X_\nu^m$ is a Markov chain on $T(\la i_1,\dots,i_{k-1}\ra )$. In
fact, until $X_.^m$ reaches $T_m$ for the first time we can view $\lbrace
X_\nu^m\rbrace$ as a realization of the Markov chain $\lbrace
X_\nu\rbrace$ of Lemma 3. 
(The equivalence classes $\la \widehat {x}\ra $ now consist
of single points $\la x\ra $ only.) Observe now that if $X_.^m$ or $X_.$ starts
from some vertex $\la i_1,\dots,i_k\ra \in T(i_1,\dots,i_{k-1})$ with
$k< s<m$, then neither $X_.^m$ nor $X_.$ can reach $T_m$ (and a fortiori
neither can reach $\infty$) without passing through some point of ${\Cal
A}_s$. Therefore
$$
\align
&\Pr\lbrace X_.^m\text{ does not reach }{\Cal A}_s\text{ between }\tau_k\text
{ and }\tau_{k-1} \mid  
\tau_k,X_0^m,X_1^m,\dots,X_{\tau_{k}}^m\rbrace\tag 2.51\\
&\hskip.5cm=\Pr\lbrace X_.\text{ reaches }\la i_1,\dots,i_{k-1}\ra 
\text{ before it reaches }{\Cal A}_s \mid 
X_0=\la i_1,\dots,i_k\ra \rbrace\\
&\hskip.5cm\le\Pr\lbrace X_.\text{ reaches }\la i_1,\dots,i_{k-1}\ra 
\text{ before it reaches }\infty\\
&\hskip5cm\text{ in }T(i_1,\dots,i_{k-1}) \mid X_0=\la i_1,\dots,i_k\ra \rbrace.
\endalign
$$
As in (2.34) the last probability equals
$$
\frac{R(T(i_1,\dots,i_k))}{R(e(i_1,\dots,i_k))+R(T(i_1,\dots,i_k))}\tag 2.52
$$
$(e(i_1,\dots,i_k)$ again denotes the edge between $\la i_1,\dots,i_{k-1}\ra $ and $\la i_1,\dots,i_k\ra )$. In particular, if there exists some $j_{k+1}\ne i_{k+1}$ such that
$$
R(T^{j_{k+1}}(i_1,\dots,i_k))\le 2L,
$$
then (2.52) is at most $2L/(2L+\varepsilon )$. Thus, if 
$$
\align
\Delta (i_1,\dots,i_s)={}&\Big\{\text{number of }k\in 
 [1,s-1]\text{ for which
there exists a }\\
& j_{k+1}\ne i_{k+1}\text{ with }
R(T^{j_{k+1}}(i_1,\dots,i_k))\le 2L\Big\},
\endalign
$$
then we obtain from the strong Markov property 
$$
\align
&\Pr\lbrace X_.^m \text{ reaches }\la 0\ra 
   \text{ before it reaches }{\Cal A}_s \mid 
 T,R,{\Cal A}_s\ne\emptyset, X_0^m=\la i_1,\dots,i_{s-1}\ra \rbrace\\
&\hskip.5cm\le\prod_{k=1}^{s-1}  \Pr\lbrace X_.\text{ does not reach }{\Cal A}_s
\text{ between }\tau_k\text{ and }\tau_{k-1} \mid 
 T,R,{\Cal A}_s\ne\emptyset,\tau_k,X_0^m,\dots,X_{\tau_{k}}^m\rbrace\\
&\hskip.5cm\le \left(\frac{2L}{2L+\varepsilon} \right)
^{\Delta (i_1,\dots,i_s)}.
\endalign
$$
In view of the definition of $\Pi$ this implies
$$
\Pi (T_{[m]},R,s)\le\bigg(\frac{2L}{2L+\varepsilon}\bigg)^{\Delta_{s}}
\text{ on }\lbrace{\Cal A}_s\ne\emptyset\rbrace ,\tag 2.53
$$
where
$$
\Delta_s=\min_{\la i_1,\dots,i_s\ra\in T_s}\Delta (i_1,\dots,i_s).
$$
Now take $L$ so large that 
$$
\Pr\lbrace R(T^j)>2L \mid \la j\ra \in T_1\rbrace\le q+2\varepsilon_0
$$
(this can be done as shown in (2.37)) and apply Lemma 1 with 
${\Cal P}$ the property of trees $t$ that $R(t)\le 2L$. By virtue of 
(2.16) there exist $0< C_1,C_2<\infty$ (which depend only on ${\Cal P}$ 
and hence only on $L$) such that 
$$
\Pr\lbrace{\Cal A}_s\ne\emptyset\text{ and }
\Delta_s\le C_1s\rbrace\le e^{-C_{2}s}\quad \text{for }s\ge 2.
$$
The lemma therefore follows from (2.53).
\qed\enddemo

\section {3. Proofs of Theorems 2 and 3}
We begin with the proof of Theorem 2. This proof is quite similar to that
of Lemma 5 in Grimmett and Kesten (1983). The result will follow fairly
easily from considerations about the random subgraphs $\tau^0$ and
$\tau^\infty$ of $K_{n+2}$ whose vertices are defined to be those which
are connected to $0$ and $\infty$ respectively, by conducting paths. We
call an edge $e$ of $K_{n+2}$ {\it conducting\/} if $R(e)<\infty$; a
{\it path\/} is called {\it conducting\/} if all its edges are
conducting. According to (1.1), $\Pr\lbrace e\text{ is conducting}\rbrace
=\gamma (n)/n$ for each $e\in K_{n+2}$, and given that $e$ is
conducting $R(e)$ has conditional distribution function $F$. The number of
conducting edges in $K_{n+2}$ incident to a fixed vertex $v$ has the
binomial distribution $B(n+1,\gamma (n)/n)$, ($(n+1)$ trials with success
probability $\gamma (n)/n$ for each trial). Now let $\tau_0^0=\lbrace
0\rbrace$ (respectively $\tau_0^\infty=\lbrace\infty\rbrace )$ and let
$\tau_k^0$ (respectively $\tau_k^\infty$) be the set of vertices of
$K_{n+2}$ which can be connected to $0$ (respectively $\infty$) by a path
of $k$ conducting edges, but not by a shorter conducting path. Clearly the
$\tau_k^0$, $k\ge 0$, are disjoint, and each vertex in $\tau_k^0$ is
connected by one conducting edge to some vertex in $\tau_{k-1}^0$. There
may be several vertices in $\tau_{k-1}^0$ for which this holds. However,
there is never a conducting edge connecting a vertex in $\tau_j^0$ with a
vertex in $\tau_k^0$ when $|k-j|\ge 2$. Similar statements hold for the
$\tau_k^\infty$. The first lemma of this section states that the graph
consisting of the vertices $\bigcup\tau_k^0$ and the conducting edges between
them converges in some distributional sense as $n\to\infty$ to a family
tree of a Bienaym\'e--Galton--Watson branching process. Let $T=T^\gamma$ be
the family tree of such a branching process, whose offspring distribution
is a Poisson distribution with mean $\gamma$ as described in the
Introduction. Also $T_n$ and $T_{[n]}$ are as described in the
Introduction. Similarly $\tau_{[k]}^0$ and $\tau_{[k]}^\infty$ are the
graphs with vertex sets $\bigcup_{m=0}^k\tau_m^0$ and
$\bigcup_{m=0}^k\tau_m^\infty$, respectively, and edge sets the sets of
conducting edges between these vertices. For a fixed rooted labeled tree
$t$ consisting of a root $\rho$ and $k$ generations, the statement
$\tau_{[k]}^0=t$ means that there exists a graph isomorphism between
$\tau_{[k]}^0$ and $t$ in which $0$ corresponds to $\rho$. A similar
definition holds for $\tau_{[k]}^\infty$ or $T_{[k]}$.

\proclaim{Lemma 7}
If (1.1) holds and $\gamma (n)\to\gamma <\infty$ then for any fixed rooted
labeled tree $t$ of $k$ generations
$$
\Pr\lbrace \tau_{[k]}^0=t\rbrace =\Pr\lbrace\tau_{[k]}^\infty =t
\rbrace\to\Pr\lbrace T_{[k]}^\gamma =t\rbrace\tag 3.1
$$ 
as $n\to\infty$. More generally, if $t_1$ and $t_2$ are two fixed rooted
labeled trees with $k$ generations then
$$
\Pr\lbrace\tau_{[k]}^0\text{\rm\ is disjoint from }\tau_{[k]}^\infty\rbrace
\to 1\text{ as }n\to \infty\tag 3.2
$$
and
$$
\Pr\lbrace\tau_{[k]}^0=t_1,\tau_{[k]}^\infty=t_2\rbrace
\to\Pr\lbrace T_{[k]}^\gamma=t_1\rbrace\Pr\lbrace 
T_{[k]}^\gamma=t_2\rbrace\tag 3.3
$$
as $n\to\infty$.
\endproclaim

\demo{Proof}
We prove
$$
\Pr\lbrace\tau_{[k]}^0=t\rbrace\to\Pr\lbrace T_{[k]}^\gamma =t\rbrace\tag 3.4
$$ 
if $\gamma (n)\to\gamma <\infty$. Since $\tau_{[k]}^0$ and
$\tau_{[k]}^\infty$ clearly have the same distribution this will prove
(3.1). It will be clear how to generalize the argument to obtain (3.2) and
(3.3).

To prove (3.4) consider for any set $A$ of vertices of $K_{n+2}$ and a
vertex $x$ of $K_{n+2}$, the number of vertices outside $A$ connected by a
conducting edge to $x$. Denote this random number by $N(x,A)$. Then, under
(1.1), $N(x,A)$ has a $B(n+2-|A\cup\lbrace x\rbrace |,\gamma (n)/n)$
distribution. Thus if $A$ varies with $n$ such that $|A_n|/n\to 0$, and
$\gamma (n)\to\gamma$, then by the familiar Poisson limit for the binomial
distribution
$$
\Pr\lbrace N(x,A_n)=\ell\rbrace\to e^{-\gamma}\frac{\gamma^\ell}{\ell
!},\quad\ell =0,1,\dots.\tag 3.5
$$
Now let $t$ be a rooted labeled tree of one generation. If the size of the
first generation of $t$ equals $\ell$, then
$$ \Pr\lbrace T_{[1]}^\gamma =t\rbrace =\Pr\lbrace   |  T_1^\gamma
|=\ell\rbrace =e^{-\gamma}\frac{\gamma^\ell}{\ell !}, 
$$
because $|T_1^\gamma |$ has a Poisson distribution with mean $\gamma$.
Also $\tau_{[1]}^0$ consists of the root $0$ and a random number of
vertices connected to $0$ by a single conducting edge. This number is
precisely $N(0,\lbrace 0\rbrace )$, so that by (3.3)
$$
\align
\Pr\lbrace\tau_{[1]}^0=t\rbrace &=\Pr\lbrace |\tau_1^0|=\ell
\rbrace =\Pr\lbrace N(0,\lbrace 0\rbrace )=\ell\rbrace\\
&\to e^\gamma\frac{\gamma^\ell}{\ell !}.
\endalign
$$
Thus (3.4) holds for $k=1$. Of course it also holds for $k=0$.

We now prove (3.4) by induction on the number of generations of $t$. Let
$t'$ be a rooted labeled tree of $(k+1)$ generations. Let $t$ be the
subtree of the first $k$ generations, and denote the vertices in the
$k$th generation by $v_1,v_2,\dots ,v_M$. Finally, let $\nu (v_i)$ be the
number of vertices in the $(k+1)$th generation of $t'$ which are
connected by an edge to $v_i$. (Thus in total $\nu (v_i)+1$ edges are
incident to $v_i$ in $t'$.) Our induction hypothesis is that (3.4) holds
for the given $k$ and $t$. To prove that (3.4) also holds when $k$ is
replaced by $k+1$ and $t$ by $t'$ we start with the trivial relation
$$
\Pr\lbrace T_{[k+1]}^\gamma =t'\rbrace =\Pr\lbrace T_{[k]}^\gamma =t
\rbrace \Pr\lbrace T_{[k+1]}^\gamma =t' \mid T_{[k]}^\gamma =t\rbrace .
$$ 
If $T_{[k]}^\gamma =t$ then there exists some isomorphism $I$ between $t$
and $T_{[k]}^\gamma$. Denote by $\la v_i\ra $ the vertex of $T_{[k]}^\gamma$
which is the image of $v_i$ under $I$. $I$ can be extended to an
isomorphism between $t'$ and $T_{[k+1]}^\gamma$ if and only if $\la v_i\ra $ has
exactly $\nu (v_i)$ children in $T_{k+1}^\gamma$, $i=1,\dots,M$. The latter
event has probability
$$
\prod_{i=1}^M e^{-\gamma}\frac{\gamma^{\nu (v_i)}}{\nu (v_i)!},\tag 3.6
$$
because the numbers of children of the $\la v_i\ra $ are independent Poisson
variables with mean $\gamma$. If $T_{[k]}^\gamma=t$ there may be several
choices for $I$, but each one can be extended to an isomorphism of
$T_{[k+1]}^\gamma$ and $t'$ if and only if each $\la v_i\ra \in
T_k^\gamma$ has the correct number of children (which depends on $I$).
Denote by $\lambda (t)$ the number of distinct assignments of children to
each $\la v_i\ra \in T_k^\gamma$ which will make $T_{[k+1]}^\gamma$
isomorphic to $t'$. The conditional probability, given $T_{[k]}^\gamma$,
of the occurrence of a specific assignment of numbers of children is given
by (3.6). This is true for each of the possible assignments, since (3.6)
depends only on $t$, and not on $I$. Consequently
$$
\Pr\lbrace T_{[k+1]}^\gamma =t' \mid T_{[k]}^\gamma =t\rbrace=
\lambda (t)\prod_{i=1}^M e^{-\gamma} \frac{\gamma^{\nu 
(v_{i})}}{\nu (v_i)!}. \tag 3.7
$$ 

We can analyze
$$
\Pr\lbrace\tau_{[k+1]}^0=t' \mid \tau_{[k]}^0=t\rbrace  \tag 3.8
$$ 
in a similar way. Let $J$ be an isomorphism between 
$t$ and $\tau_{[k]}^0$ and let $x_i$ be the image in $\tau_{k}^0$ 
of $v_i$. Denote by $A$ the set of vertices of $K_{n+2}$ which 
belong to $\tau_{[k]}^0$. Then $J$ can be extended to an 
isomorphism of $t'$ and $\tau_{[k+1]}^0$ if and only if
$$
\align
&x_i\text{ is connected by a conducting edge to exactly }
\nu (v_i)\text{ vertices }\tag "\rlap{(3.9)}"\\
&\text{of }K_{n+2}\text{ outside }A,\text{ but no }
x\in K_{n+2}\setminus A\text{ is connected by conducting }\\
&\text{edges to
 two of the vertices }x_1,\dots,x_M.
\endalign
$$
Denote by $B_i$ the collection of vertices outside $A$ which are connected
by a conducting edge to $x_i$. Then (3.9) is just the event
$$
\lbrace |B_i|=\nu (v_i),1\le i\le M,\text{ and the }B_i\text{ are disjoint}
\rbrace.
$$  
Thus, if we set $A_0=A$, $A_i=A_0\cup B_1\cup\dots\cup B_i$, then the
conditional probability of (3.9) given $\tau_{[k]}^0=t$ can be written as
$$
\align
\prod_{i=1}^M \Pr\Big\lbrace  |B_i|=\nu (v_i),{}&B_i\subset K_{n+2}\setminus
A_{i-1} \,\Big|\, \tau_{[k]}^0=t, \tag 3.10\\ 
&|B_j|=\nu (v_j),1\le j\le i-1\text{ and
}B_1,\dots,B_{i-1}\text{ are disjoint}\Big\rbrace . 
\endalign
$$
Finally, knowledge of $\tau_{[k]}^0$ and $B_1,\dots,B_{i-1}$ gives no
information about edges between $x_i$ and $K_{n+2}\setminus A$. Therefore
the $i$th factor in (3.10) equals
$$
\binom{n+2-|A_{i-1}|}{\nu(v_i)} \bigg(\frac{\gamma 
(n)}{n}\bigg)^{\nu (v_i)}\bigg(1-\frac{\gamma (n)}{n}\bigg)^{n+2-|A|-\nu 
(v_i)}. \tag 3.11
$$
As in (3.5) the limit of (3.11) as $n\to\infty$ equals
$$
e^{-\gamma}\frac{\gamma^{\nu (v_{i})}}{\nu (v_i)!}
$$
so that (3.10) converges to (3.6). It follows that (3.8) converges to the
right hand side of (3.7). This, together with the induction hypothesis
implies
$$
\Pr\lbrace\tau_{[k+1]}^0=t'\rbrace\to\Pr\lbrace T_{[k+1]}^\gamma =t'\rbrace .
$$ 
This completes the induction step for the proof of (3.4). As mentioned
before, the proofs of (3.2) and (3.3) follow along similar lines when
$\gamma (n)\to\gamma <\infty$.
\qed\enddemo

\proclaim{Lemma 8}
Assume (1.1). If $\gamma (n)\to \gamma <\infty$, then 
$$
\limsup\limits_{n\to\infty} \Pr\lbrace R_n\le x\rbrace\le \Pr\lbrace
R'(\gamma )+R''(\gamma )\le x\rbrace, \tag 3.12
$$
at each continuity point of the right hand side, where $R'(\gamma )$ and
$R''(\gamma )$ are independent random variables, each with the
distribution of $R(T^\gamma )$. 
\endproclaim

\demo{Proof}
Assume that for some $n$ and a certain realization of the resistances and
fixed $k$, $\tau_{[k]}^0$ and $\tau_{[k]}^\infty$ are disjoint. By the
monotonicity property (2.6), $R_n$, the resistance between $0$ and
$\infty$ in $K_{n+2}$, is then at least equal to the resistance between
$0$ and $\infty$ in the network consisting of
$\tau_{[k]}^0\cup\tau_{[k]}^\infty$ and shortcircuits between all pairs of
vertices $v,w\in
K_{n+2}\setminus\tau_{[k-1]}^0\cup\tau_{[k-1]}^\infty$. Indeed the
resistance of each edge in this network is less than or equal to the
resistance assigned to it in the original network on $K_{n+2}$, since by
the definition of $\tau^0$ and $\tau^\infty$ there are no conducting edges
in the original network between $\tau_{[k-1]}^0\cup\tau_{[k-1]}^\infty$
and $K_{n+2}\setminus\tau_{[k-1]}^0\cup\tau_{[k-1]}^\infty$, except for
the edges between $\tau_{k-1}^0$ and $\tau_k^0$ and between
$\tau_{k-1}^\infty$ and $\tau_k^\infty$.

\topinsert
\figure
\mletter{0}{1.5}{1.9}
\mletter{\tau^0_{[3]}}{3}{3.1}
\mletter{\infty}{10.8}{1.9}
\lastletter{\tau^\infty_{[3]}}{10}{3.1}
\centerline{\epsfxsize=9cm\epsffile{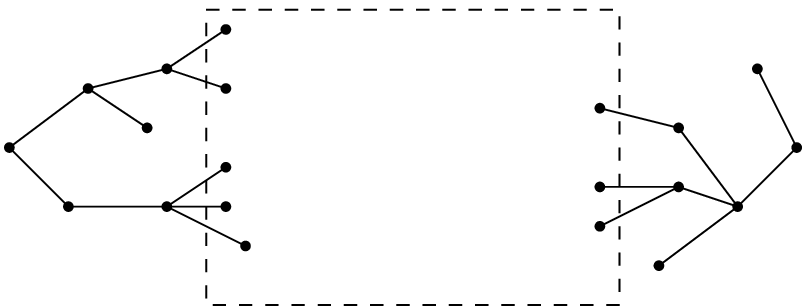}}

\capt{7}{All vertices in the dashed rectangle are
shortcircuited; this illustrates the 
construction in the proof of Lemma 8, with $k=3$.}
\endinsert

Now the network with the shortcircuits inserted is equivalent to
$\tau_{[k]}^0$ and $\tau_{[k]}^\infty$ in series, after all vertices in
$\tau_k^0\cup\tau_k^\infty$ are identified as a single vertex (see Figure
7 for $k=3$ and Figure 8).

\topinsert
\figure
\mletter{0}{1.5}{1.3}
\lastletter{\infty}{10.8}{1.3}
\centerline{\epsfxsize=9cm\epsffile{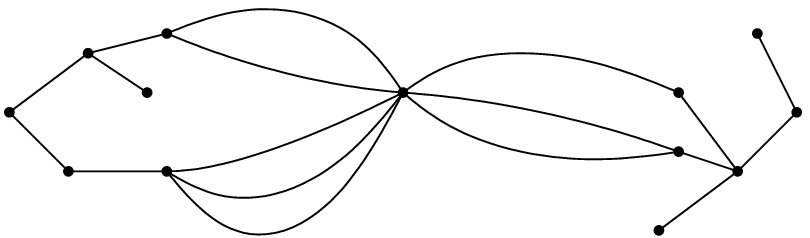}}

\capt{8}{The network of Figure 7 is equivalent
to the above network.}
\endinsert

Denote by $r_k^0$ the resistance between $0$ and $\tau_k^0$ in
$\tau_{[k]}^0$ when all vertices in $\tau_k^0$ are identified (or
shortcircuited). Define $r_k^\infty$ similarly by replacing $0$ by
$\infty$. Then the resistance between $0$ and $\infty$ in the network with
short circuits is $r_k^0+r_k^\infty$. By Lemma 7 the probability that
$\tau_{[k]}^0$ and $\tau_{[k]}^\infty$ are disjoint tends to $1$, while
$\tau_{[k]}^0$ and $\tau_{[k]}^\infty$ converge in distribution to two
independent trees with the distribution of $T_{[k]}^\gamma$. Moreover, the
fact that an edge belongs to $\tau_{[k]}^0\cup\tau_{[k]}^\infty$ says no
more about its resistance than that this resistance is finite. Thus the
resistance of each of the edges of $\tau_{[k]}^0\cup\tau_{[k]}^\infty$ has
the distribution function $F$. In addition these resistances are
independent. Therefore $(r_k^0,r_k^\infty )$ converges in distribution to
$(R_k'(\gamma )),R_k''(\gamma ))$, where $R_k',R_k''$ are independent,
each with the distribution of
$$
R(T_{[k]}^\gamma )=\{\text{resistance between }\la 0\ra \text{ and }T_k^\gamma
\text{ in }T_{[k]}^\gamma\} .
$$    
In view of the above
$$
\limsup\Pr\lbrace R_n\le x\rbrace\le\limsup \Pr\lbrace
r_k^0+r_k^\infty\le x\rbrace
=\Pr\lbrace R_k'(\gamma )+R_k''(\gamma
)\le x\rbrace ,
$$
at each continuity point of the last member. But $R(T^\gamma
)=\lim\limits_{k\to\infty }R(T_{[k]}^\gamma )$ by definition, so that
(3.12) follows.
\qed\enddemo

Note that Lemma 8 proves ``one half" of Theorem 3. It shows that $R_n$ is
asymptotically at least as large in distribution as $R'(\gamma )+R''(\gamma
)$. Also, a much simplified version of the proof of Lemma 8 implies
Theorem 2, as we now show.

\demo{Proof of Theorem 2} $R_n=\infty$ if $0$ is not connected to $\infty$
by a conducting path. Thus, by Lemma 7,
$$
\align
\liminf\limits_{n\to\infty}\Pr\lbrace
R_n=\infty\rbrace &\ge\liminf\limits_{n\to\infty} \Pr\lbrace\text{for some
}k<\infty,\tau_k^0=\emptyset\text{ and
}\infty\notin\tau_{[k]}^0\rbrace\\
&\ge \Pr\lbrace T_k^\gamma
=\emptyset\text{ for some }k\rbrace.
\endalign
$$
However, the last probability is just the extinction probability of the
branching process $\lbrace Z_n\rbrace$, and this probability equals $1$ if
$\gamma\le 1$ (see Harris (1963) Theorem I.6.1). 
\qed\enddemo

To obtain an upper bound for $R_n$ in Theorem 3 we need to show that
$\tau_{[k]}^0$ is close in distribution to $T_{[k]}^\gamma$ not only for
fixed $k$ and large $n$, but even for $k$ a suitable multiple of $\log n$.
Since we only want an upper bound for $R_n$, it suffices (as we shall see)
to show that for $k$ a suitable multiple of $\log n$ and for fixed $\delta
<\gamma$, $\tau_{[k]}^0$ is stochastically larger than $T_{[k]}^\delta$. We
shall do this by ``coupling". For the remainder of these notes we assume
(1.4) and take $1<\delta <\gamma$. We shall construct on one probability
space $\tau_{[k]}^0$ and two other graphs $\Tilde{\tau}_{[k]}^0$ and
$\Tilde{T}_{[k]}=\Tilde{T}_{[k]}^\delta$, such that
$\Tilde{T}_{[k]}$ has the same distribution as $T_{[k]}^\delta$ and
such that $\Tilde{\tau}_{[k]}^0$ and $\Tilde{T}_{[k]}$ are trees with root
at $0$ and such that with high probability
$$
\Tilde{T}_{[k]}\subset\Tilde{\tau}_{[k]}^0\subset\tau_{[k]}^0. \tag 3.13
$$ 
$\tau_{[k]}^0$ has already been constructed in the beginning of this
section. We construct $\Tilde{\tau}_{[\cdot]}^0$ as a subgraph of
$K_{n+2}\setminus\lbrace\infty\rbrace$ in stages. We set $B_0=\lbrace
0\rbrace =\Tilde{\tau}_0^0$. $\Tilde{\tau}_{[0]}^0$ is the graph which
consists of the vertex $0$ only. At stage $\ell$, $\Tilde{\tau}_{[j]}^0$ will
have been constructed for $j\le\ell$, such that these graphs are an
increasing family of trees in $K_{n+2}\setminus\lbrace\infty\rbrace$. Let
$B_j$ denote the set of vertices of $\Tilde{\tau}_{[j]}^0$, and let
$i_1<\dots <i_r$ be the vertices in $B_\ell\setminus B_{\ell -1}$. (If
$\ell\ge 1$, then $i_1,\dots,i_r\in\lbrace 1,\dots ,n\rbrace.)$ We
then construct $\Tilde{\tau}_{[\ell +1]}^0$ by choosing disjoint subsets
of $\lbrace 1,2,\dots ,n\rbrace\setminus B_\ell$ for the vertices which
will be connected by an edge in $\Tilde{\tau}_{[\ell +1]}^0$ to $i_1,\dots
,i_r$. These choices too are made successively. Let $C_\ell(p)$ denote the
union of $B_\ell$ and all vertices of $\lbrace 1,\dots,n\rbrace$ which
have been chosen already to be connected to $i_1,i_2,\dots,i_{p-1}$. Thus
$C_\ell(1)=B_\ell$. We now add to the vertex set of $\Tilde{\tau}_{[\ell
+1]}^0$ all vertices in $\lbrace 1,\dots,n\rbrace\setminus C_\ell(p)$
which are connected by a {\it conducting\/} edge to $i_p$. The edges
between these vertices and $i_p$ are added to the edge set of
$\Tilde{\tau}_{[\ell +1]}^0$. $\Tilde{\tau}_{[\ell +1]}^0$ is the graph
obtained after all these additions have been performed for
$i_1,\dots,i_r$. It is clear from the construction that
$\Tilde{\tau}_{[\ell ]}^0$ is a tree for each $\ell$, since each time we
only add vertices which have not been used before and only one edge
between each new vertex and the old vertices. Comparison with the
construction of $\tau_{[\ell ]}^0$ also shows immediately that
$\Tilde{\tau}_{[\ell ]}^0$ is a subgraph of $\tau_{[\ell ]}^0$.

To construct $\Tilde{T}_{[\cdot ]}$ it is convenient to view $K_{n+2}$ as
a subgraph of $K_\infty$, which is the complete graph with vertices
$0,\infty,1,2,\dots$. (Recall that $K_{n+2}$ has vertices
$0,\infty,1,2,\dots,n$.) We shall use an auxiliary family of random
variables
$\lbrace Y(j),U(j):j\ge 0\rbrace$, with each $Y$ a Poisson
variable with mean $\delta$ ($1<\delta <\gamma$ is a fixed number) and each
$U$ uniformly distributed on $[0,1]$. All these variables are taken
independent of each other and independent of all resistances and of all
$\Tilde{\tau}_{[\ell ]}^0$. Again we construct $\Tilde{T}_{[k]}$ in
stages. Set $D_0=\lbrace 0\rbrace=\Tilde{T}_0$ and take for
$\Tilde{T}_{[0]}$ the graph consisting of the vertex $0$ only. Assume we
have already chosen $\Tilde{T}_{[j]}$, $j\le\ell$ as subgraphs of $K_\infty$
such that for each $j\le\ell$ each vertex of $\Tilde{T}_{[j]}$ either is a
vertex of $\Tilde{\tau}_{[j]}^0$, or belongs to $\lbrace
n+1,n+2,\dots\rbrace\subset K_\infty\setminus K_{n+2}$. We then choose
$\Tilde{T}_{\ell +1}$ as follows. Let $D_\ell$ be the vertex set of
$\Tilde{T}_{[\ell ]}$, and let $j_1<j_2<\dots <j_s$ be the vertices in
$D_\ell\setminus D_{\ell -1}$. Again we add successively disjoint sets of
vertices and connect them by edges to $j_1,j_2,\dots ,j_s$, respectively,
to form $\Tilde{T}_{[\ell +1]}$. Denote by $E_\ell(p)$ the union of
$D_\ell$ and all the vertices which have already been connected to
$j_1,\dots ,j_{p-1}$; $E_\ell (1)=D_\ell$. We now choose the vertices
connected to $j_p$. First consider the case where $j_p\le n$. By our
inductive assumption $j_p$ is then a vertex of $\Tilde{\tau}_{[\ell ]}^0$,
since it belongs to $\Tilde{T}_{[\ell ]}$ as well as to $\lbrace 1,\dots
,n\rbrace$. Let $j_p=i_\nu$ and let $r_1<r_2<\dots <r_q$ be the vertices of
$\Tilde{\tau}_{[\ell +1]}^0$ which are connected by an edge of
$\Tilde{\tau}_{[\ell +1]}^0$ to $i_\nu =j_p$. By construction all
$r_i\in\lbrace 1,\dots,n\rbrace\setminus C_\ell (\nu)$. Note that
$q=0$ is possible, so that there may not be any vertices of this kind. Put
$\beta (-1)=0$ and for $x\ge 0$
$$
\align
\beta(x)&=\beta(x;n,|C_\ell(\nu)|)\tag "\rlap{(3.14)}"\\ 
&=\sum_{j\le x} \binom{n+1-
|C_{\ell}(\nu )|}{j} \bigg(\frac{\gamma (n)}{n}\bigg)^j
\bigg(1-\frac{\gamma(n)}{n}\bigg)^{n+1-|C_\ell (\nu)|-j}
\endalign  
$$
and
$$
\pi(x)=\pi(x;\delta)=\sum\limits_{j\le x}e^{-\delta}\frac{\delta^{j}}{j!}.
$$
With $U$ the previously chosen uniform random variable we take 
$$
\align
V(j_p)&=\beta(q-1)+U(j_p)[\beta(q)-\beta(q-1)], \tag 3.15\\
u&=\pi^{-1}(V(j_p)):= \{\text{smallest integer }m\text{ for which }
\pi(m)\ge V(j_p)\}. 
\endalign
$$
If $u\le q$, then we add to $\Tilde{T}_{[\ell +1]}$ the vertices
$r_1,\dots,r_u$ plus the edges between $r_1,\dots,r_u$ and $j_p=i_\nu$
(all these edges are also edges of $\Tilde{\tau}_{[\ell+1]}^0$). If $u>q$,
then we add $r_1,\dots,r_q$ plus the edges between these and $j_p$, and in
addition choose $u-q$ vertices from $\lbrace
n+1,n+2,\dots\rbrace\setminus E_\ell(p)$ and add these vertices also to
$\Tilde{T}_{[\ell +1]}$, together with an edge from each of them to $j_p$.
Thus in each case $u$ new vertices are connected to $j_p$ in
$\Tilde{T}_{[\ell +1]}$. Finally, if $j_p>n$, then we add $Y(j_p)$
vertices from $\lbrace n+1,n+2,\dots\rbrace\setminus E_\ell (p)$ to
$\Tilde{T}_{[\ell+1]}$ and an edge from each of these vertices to
$j_p$ ($Y(\cdot )$ is the previously chosen Poisson variable).
$\Tilde{T}_{[\ell +1]}$ is the graph obtained when the above construction
is completed for all $j_1,\dots,j_s.$

\proclaim{Lemma 9}
Let $1<\delta <\gamma$ be fixed and let  
\footnote{$\lfloor a\rfloor$ denotes the largest integer $\le a$.} 
$$
m=m_n=\left\lfloor\frac{3}{4}\frac{\log n}{\log\gamma}\right\rfloor.\tag 3.16
$$
Then $\Tilde{T}_{[m]}^\delta$ is a tree with the same distribution as
$T_{[m]}^\delta$. Also $\Tilde{\tau}_{[m]}^0$ is a tree, and as
$n\to\infty$,
$$
\Pr\lbrace\Tilde{T}_{[m]}^\delta\subset\Tilde{\tau}_{[m]}^0\subset\tau_{[m]}^0\rbrace\to
1.\tag 3.17
$$
\endproclaim

\demo{Proof}
To show that $\Tilde{T}_{[m]}^\delta$ is a tree we merely have to observe
that at each stage of its construction we add vertices which have not been
used before and one edge from each new vertex to one old vertex. Thus at
no stage can a circuit arise in any $\Tilde{T}_{[\ell ]}$.

To show that $\Tilde{T}_{[m]}$ has the same distribution as
$T_{[m]}^\delta$ we must show that for $\ell<m$ the ``number of children"
of each vertex of $\Tilde{T}_{[\ell]}$ in $\Tilde{T}_{[\ell+1]}$ has a
Poisson distribution with mean $\delta$, and that all these numbers are
independent. Use the notation of the construction preceding this lemma.
Let $j_p$ be a vertex in $\Tilde{T}_{[\ell]}$. If $j_p>n$ then its
children (i.e., vertices of
$\Tilde{T}_{[\ell+1]}\setminus\Tilde{T}_{[\ell]}$ connected to $j_p$) are
precisely $Y(j_p)$ vertices from $\lbrace n+1,n+2,\dots\rbrace\setminus
E_\ell (p)$. Since $Y(j_p)$ was a Poisson variable with mean $\delta$
independent of all other $Y$ and $U$, there is nothing to prove in this
case. Now consider a $j_p\in\lbrace 1,\dots ,n\rbrace$ with
$j_p=i_\nu$ as above. In this case $j_p$ has $u=\pi^{-1}(V(j_p))$
children. The distribution of $u$ is given by
$$
\Pr\lbrace u\le r \mid \Tilde{\tau}_{[j]}^0,\Tilde{T}_{[j]},j\le\ell,C_\ell
(\nu )\rbrace
=\Pr\lbrace V(j_p)\le\pi
(r) \mid \Tilde{\tau}_{[j]}^0,\Tilde{T}_{[j]},j\le\ell,C_\ell (\nu)\rbrace.\tag
3.18
$$
Recall the definition of $V$ in (3.15) and note that $q$ in this formula
is just the number of conducting edges between $j_p=i_\nu$ and vertices in
$\lbrace 1,\dots,n\rbrace\setminus C_\ell (\nu)$. Given
$\Tilde{\tau}_{[j]}^0$, $\Tilde{T}_{[j]}$, $j\le\ell $, and the sets
$C_\ell(\nu)$, $E_\ell(p)$, the conditional distribution of $q$ is binomial
$B(n+1-|C_\ell (\nu)|,\gamma (n)/n)$. In particular, if $0<\gamma (n)<n$
and $q_0<n$ is a fixed integer, then
$$
\Pr\lbrace\beta (q)\le\beta
(q_0)\mid \Tilde{\tau}_{[j]}^0,\Tilde{T}_{[j]},j\le\ell,C_\ell (\nu)\rbrace
=\Pr\lbrace q\le q_0 \mid \Tilde{\tau}_{[j]}^0,\Tilde{T}_{[j]},j\le\ell,C_\ell
(\nu)\rbrace.\tag3.19
$$
If we take
$$
q_0=q_0(r)=\lbrace\text{largest integer }s\text{ with }\beta(s)\le\pi(r)\rbrace
$$
then we obtain from (3.15), (3.18) and (3.19) 
$$
\align
&\Pr\lbrace u\le
r \mid \Tilde{\tau}_{[j]}^0,\Tilde{T}_{[j]},j\le\ell,C_\ell (\nu)\rbrace\\
&\hskip1cm=\Pr\lbrace q\le q_0 \mid \Tilde{\tau}_{[j]}^0,\Tilde{T}_{[j]},j\le\ell, C_\ell
(\nu)\rbrace\\ 
&\hskip1.5cm +\Pr\lbrace
q=q_0+1 \mid \Tilde{\tau}_{[j]}^0,\Tilde{T}_{[j]},j\le\ell,C_\ell(\nu)\rbrace
\Pr\left\lbrace U(j_p)\le \frac{\pi (r)-\beta (q_0)}{\beta
(q_{0}+1)-\beta (q_0)}\right\rbrace\\ 
&\hskip1cm=\beta
(q_0)+\pi(r)-\beta(q_0)=\pi(r).
\endalign
$$
Thus the number of children of any vertex in $\Tilde{T}_{[\ell ]}$ indeed
has a Poisson distribution with mean $\delta$, as desired. A slightly
closer look at the above argument shows also that the numbers of children
of each of the vertices in $\Tilde{T}_{[\ell ]}$ are independent, so that
the first claim of the lemma follows.

\indent We already observed that $\Tilde{\tau}_{[m]}^0$ is a subtree of
$\tau_{[m]}^0$ by construction, so that we only need to prove
$$
\Pr\lbrace\Tilde{T}_{[m]}\subset\Tilde{\tau}_{[m]}^0\rbrace\to 1 \tag 3.20
$$
for (3.17). First observe that
$\Tilde{T}_{[m]}\subset\Tilde{\tau}_{[m]}^0$ fails only if there exists
some $\ell<m$ and a vertex $j_p$ of $\Tilde{T}_{[\ell]}$ which equals a
vertex $i_\nu$ of $\Tilde{\tau}_{[\ell]}^0$ such that $u>q$ where $u$ is
the number of children of $j_p$ in $\tilde{T}_{[\ell +1]}$ and $q$ is the
number of children of $i_\nu$ in $\Tilde{\tau}_{[\ell +1]}^0$. By our
construction this requires $\pi (q;\delta)<\beta (q;n,|C_\ell (\nu )|)$;
see (3.15). Next we obtain a lower bound for $\beta (q)$. The expected
number of vertices in any subset of $\lbrace 1,\dots,n\rbrace$ connected
by a conducting edge to any fixed vertex of $K_{n+2}$ is at most
$n\cdot\gamma (n)/n=\gamma (n)$. It follows from this that
$$
\EE\lbrace |\Tilde{\tau}_{\ell
+1}^0|\mid \Tilde{\tau}_j^0,j\le\ell\rbrace\le\gamma (n)|\Tilde{\tau}_\ell^0|
$$
and (see (3.16))
$$
\EE|\Tilde{\tau}_{[m]}^0|\le\sum\limits_{\ell=0}^m\lbrace\gamma
(n)\rbrace^\ell\le\frac{\gamma (n)}{\gamma (n)-1}n^{3/4}. $$
Therefore
$$ \Pr\lbrace |\Tilde{\tau}_{[m]}^0|>n^{7/8}\rbrace\le\frac{\gamma
(n)}{\gamma (n)-1} n^{-1/8}\to 0.\tag 3.21 
$$
If $|\Tilde{\tau}_{[m]}^0|\le n^{7/8}$, then also $|C_\ell (\nu )|\le
|\Tilde{\tau}_{[m]}^0|\le n^{7/8}$ for all $C_\ell (\nu )$ used in the
construction of $\Tilde{\tau}_{[m]}^0$. Therefore, if we set
$$
\Tilde{\gamma}(n)=\frac{n+1-n^{7/8}}{n}\gamma (n),
$$
then for any $x\le n^{1/16}$ and some constant C
$$
\align
\beta (x)&=\beta (x;n,|C_\ell (\nu)|)\\
&\le\sum_{j\le x}\binom{n+1-n^{7/8}}{j} \bigg(\frac{\gamma
(n)}{n}\bigg)^j \bigg(1-\frac{\gamma (n)}{n}\bigg)^{n+1-n^{7/8}-j}\\ &\le
(1+Cn^{-15/16})\sum\limits_{j\le
x}e^{-\Tilde{\gamma}(n)}\frac{\{\Tilde{\gamma}(n)\}^j}{j!}\\ &\le\pi
(x;\Tilde{\gamma}(n))+Cn^{-15/16}=\int_{\Tilde{\gamma}(n)}^\infty
e^{-z}\frac{z^x}{x!}\,dz+Cn^{-15/16}.
\endalign
$$
Moreover, for $x\ge\Tilde{\gamma}(n)$
$$
\int\limits_{\Tilde{\gamma}(n)}^\infty
e^{-z}\frac{z^x}{x!}\,dz+Cn^{-15/16}\le\int\limits_\delta^\infty
e^{-z}\frac{z^x}{x!}\,dz=\pi(x;\delta) 
$$ 
as long as 
$$
\int\limits_\delta^{\Tilde{\gamma}(n)} e^{-z}\frac{z^x}{x!}\,dz\ge
(\Tilde{\gamma}(n)-\delta)e^{-\delta}\frac{\delta^x}{x!}\ge Cn^{-15/16}.
$$
Therefore, if we define 
$$
 s(n)=\left\{\text{smallest }s\text{ with
}e^{-\delta}\frac{\delta^s}{s!}< \frac{2Cn^{-15/16}}{\gamma-\delta}\right\},
\tag
3.22 
$$
then for sufficiently large $n$, $\pi(q;\delta)<\beta (q;
n, |C_\ell (\nu)|)$ can occur only if \footnote{$a\wedge
b=\min\lbrace a,b\rbrace $} $q\ge s(n)\wedge n^{1/16}$. It follows that
$\Tilde{T}_{[m]}\subset\Tilde{\tau}_{[m]}^0$ whenever
$|\Tilde{\tau}_{[m]}^0|\le n^{7/8}$ and all vertices in
$\Tilde{T}_{[m-1]}$ have fewer than $s(n)\wedge n^{1/16}$ children. Thus,
by virtue of (3.21) 
$$
\align
&\Pr\lbrace \Tilde{T}_{[m]}\text{ is
not a subgraph of }\Tilde{\tau}_{[m]}^0\rbrace\tag 3.23\\ 
&\hskip1cm\le\frac{\gamma (n)}{\gamma (n)-1}n^{-1/8} \\
&\hskip3cm+\Pr\lbrace\text{some vertex of
}\Tilde{T}_{[m-1]}\text{ has 
at least }s(n)\wedge n^{1/16}\text{ children}\rbrace.
\endalign 
$$ 

Finally we use the fact that $\Tilde{T}_{[m]}$ has the same distribution
as $T_{[m]}^\delta$, so that from standard branching process formulae (see
Harris (1963) Theorem I.5.1) 
$$
\Pr\lbrace |T_{[m]}^\delta|\ge
A\delta^m\rbrace\le\frac{1}{A\delta^m} \EE|T_{[m]}^\delta
|\le\frac{\delta}{A(\delta -1)}. 
$$
 For the right hand side of (3.23) we
therefore find for large $n$ the estimate 
$$
\align
 &\frac{\gamma
(n)}{\gamma (n)-1}n^{-1/8}+\Pr\lbrace |T_{[m]}^\delta |\ge
A\delta^m\rbrace\tag 3.24\\ 
&\hskip.5cm+\Pr\lbrace\text{one of }A\delta^m\text{ independent
Poisson variables, mean $\delta$,  
 is at least
}s(n)\wedge n^{1/16}\rbrace\\ 
&\hskip1cm\le\frac{\gamma (n)}{\gamma
(n)-1}n^{-1/8}+\frac{\delta}{A(\delta -1)}+A\delta^m\sum_{k\ge
s(n)\wedge n^{1/16}} e^{-\delta}\frac{\delta^k}{k!}.
\endalign 
$$
Finally, for $t\ge 2\delta$ 
$$
\sum_{k\ge t}
e^{-\delta}\frac{\delta^k}{k!}\le \left(1-\frac{\delta}{t}\right)^{-1}
e^{-\delta}\frac{\delta^t}{t!}\le  2e^{-\delta}\frac{\delta^t}{t!} 
$$
so that by virtue of (3.22) and (3.16) 
$$
\align 
A\delta^m
\sum_{k\ge s(n)} e^{-\delta}\frac{\delta^k}{k!} &\le 2A\delta^m
e^{-\delta}\frac{\delta^{s(n)}}{s(n)!}\\ 
&\le 2A\delta^m
\frac{2C}{\gamma -\delta} n^{-15/16}\\
&=\O(n^{3/4-15/16})=\O(n^{-3/16}).
\endalign 
$$
Obviously, also 
$$
 A\delta^m\sum_{k\ge n^{1/16}}
e^{-\delta}\frac{\delta^k}{k!}\to 0\quad\text{ as }n\to\infty. 
$$
Thus the
right hand side of (3.24) can be made as small as desired by choosing
first $A$ and then $n$ large. (3.20) and (3.17) follow. 
\qed\enddemo

In the same way as we constructed $\Tilde{\tau}_{[m]}^0$ we can
construct a subtree $\Tilde{\tau}_{[m]}^\infty$ of $\tau_{[m]}^\infty$. We
want $\Tilde{\tau}_{[m]}^\infty$ disjoint from $\Tilde{\tau}_{[m]}^0$.
This can be achieved by first constructing $\Tilde{\tau}_{[m]}^0$ and then
choosing for $\Tilde{\tau}_{[m]}^\infty$ only vertices of
$K_{n+2}\setminus\Tilde{\tau}_{[m]}^0$. We can then construct a random
tree $\wwtilde T_{[m]}$ which has the same relation to
$\Tilde{\tau}_{[m]}^\infty$ as $\Tilde{T}_{[m]}$ to
$\Tilde{\tau}_{[m]}^0$. We shall, however, use Poisson and uniform
variables for $\wwtilde T$ which are independent of the $Y(\cdot
)$ and $U(\cdot)$ used in the construction of $\Tilde{T}_{[m]}$. Also
we shall choose all vertices of $\wwtilde T_{[m]}$ disjoint from
those of $\Tilde{T}_{[m]}$. It is then not hard to show that
$\Tilde{T}_{[m]}$ and $\wwtilde T_{[m]}$ are disjoint trees, which
are independent of each other, each with the distribution of
$T_{[m]}^\delta$. Moreover,
$$
\Pr\lbrace T_{[m]}^\delta\subset
\Tilde{\tau}_{[m]}^0\subset\tau_{[m]}^0\text{ and }
\wwtilde T_{[m]}^\delta
\subset\Tilde{\tau}_{[m]}^\infty\subset\tau_{[m]}^\infty\rbrace\to
1\quad\text{as } n\to\infty.\tag 3.25
$$

The next lemma is almost immediate from (3.25), but we need some more
notation. Let $T'$ and $T''$ be independent disjoint trees, each with the
distribution of $T^\delta$. $T'_m,T'_{[m]},T''_{m}$ and $T''_{[m]}$ have
the obvious meaning. The vertices of $T'$ and $T''$ will be labeled
$\la 0\ra ',\la i_1,\dots,i_n\ra '$ and 
$\la 0\ra '',\la i_1,\dots,i_n\ra ''$ in the usual way.
All edges of $T'$ and $T''$ are assigned a random resistance, chosen
according to the distribution $F$. All these resistances are assumed
independent. For each $k$ we can form a network $N(n,k)$ consisting of
$T'_{[k]},T''_{[k]}$ and a resistance between each pair of vertices
$v',v''$ with $v'\in T'_k$, $v''\in T''_k$. The latter resistances are chosen
independent of each other and of $T',T''$ and the resistances in $T'\cup
T''$. Each of the resistances between $T_{k}'$ and $T_{k}''$ is chosen
according to (1.1). Of course this is equivalent to connecting each given
pair $v',v''$ only with probability $\gamma(n)/n$, but with $F$ for the
conditional distribution function of the resistance between them, when it
is given that they are connected. We put
$$
\rho(n,k)=\{\text{resistance between }
\la 0\ra '\text{ and }\la 0\ra ''\text{ in }N(n,k)\}.
$$

\proclaim{Lemma 10} 
Let $1<\delta <\gamma$ and let $m=m_n$ be as in (3.16). Then
$$
\Pr\lbrace R_n\le x\rbrace\ge\Pr\lbrace\rho(n,m_n)\le x\rbrace+\o(1),\tag 3.26
$$
where $\o(1)\to 0$ as $n\to\infty$.
\endproclaim

\demo{Proof}
The construction in stages of $\Tilde{\tau}^0_{[m]}$, $\Tilde{T}_{[m]}$,
$\Tilde{\tau}^{\infty}_{[m]}$ and $\wwtilde T_{[m]}$ is such that it
gives no information about the resistances of edges of $K_{n+2}$ between the
last generation of $\Tilde{\tau}^0_{[m]}$ (i.e.,
$\Tilde{\tau}^0_{[m]}\setminus\Tilde{\tau}^0_{[m-1]})$ and the last
generation of $\Tilde{\tau}^\infty_{[m]}$ (i.e.,
$\tau^\infty_{[m]}\setminus\tau^\infty_{[m-1]}$). Consequently,
conditionally on
$\Tilde{\tau}^0_{[m]}$, $\Tilde{\tau}^\infty_{[m]}$,
$\Tilde{T}_{[m]}$, $\wwtilde
T_{[m]}$, all these edges have independent resistances, each with
distribution given in (1.1). We also insert an edge between any pair of
vertices $v, w$ with
$v\in\Tilde{T}_m:=\Tilde{T}_{[m]}\setminus\Tilde{T}_{[m-1]}$ and
$w\in\wwtilde T_m:=\wwtilde T_{[m]}\setminus \wwtilde
T_{[m-1]}$ and not both $v$ and $w$ in $K_{n+2}$. These edges are also
given resistances with the distribution (1.1), and we take all these
resistances independent of each other and of the ones in $K_{n+2}$.
Finally we choose further independent resistances with distribution
function $F$ for all edges in $\Tilde{T}_{[m]}\cup\wwtilde T_{[m]}$
which are not edges of $K_{n+2}$. Recall that all edges of
$\Tilde{\tau}^0_{[ m]}$ and $\Tilde{\tau}^\infty_{[m]}$ are conducting by
construction, and hence have conditional distribution $F$ for their
resistance. Thus, the resistance, $\Tilde{\rho}$ say, between $0$ and
$\infty$ in the network consisting of $\Tilde{T}_{[m]}$, $\wwtilde
T_{[m]}$ and the edges between $\Tilde{T}_{m}$ and $\wwtilde
T_{m}$ has precisely the distribution of $\rho(n,m)$. Moreover, when
$\Tilde{T}_{[m]}$ and $\wwtilde T_{[m]}$ are subgraphs of
$\Tilde{\tau}^0_{[m]}$ and $\Tilde{\tau}^ \infty_{[m]}$, respectively,
then this network is part of $K_{n+2}$. Since the resistance between $0$
and $\infty$ in any sub-network of $K_{n+2}$ is at least $R_n$ (by the
monotonicity property (2.6)) we have $\Pr\lbrace\Tilde{\rho}\ge
R_n\rbrace\to 1$ (by (3.25)). (3.26) follows because $\Tilde{\rho}$ and
$\rho(n,m)$ have the same distribution. 
\qed\enddemo

We shall now show that $\rho(n,m_n)$ converges in distribution to
$R'(\delta)+R''(\delta)$. Except for the proof of (3.72) which occurs
almost at the end of these notes we make no further use of the fact that
the offspring distribution in the branching process is a Poisson
distribution.

\proclaim{Lemma 11}
Assume that
$$
F(K)=1 \quad \text{for some }\quad K<\infty.\tag 3.27
$$
Then there exist constants $0<C_3,C_4<\infty$ such that for
$0\le\varepsilon<\frac{1}{3}$ and $k\ge\frac{1+2\varepsilon}{2\log
\delta}\log n$ and $n$ sufficiently large
$$
\align
&\Pr\lbrace\rho(n,k)\le(2k+1)K\mid |T'_k|\ne 0, |T''_k|\ne 0\rbrace\tag 3.28\\
&\hskip1cm=\Pr\lbrace\exists\text{\rm\ conducting path between $\la 0\ra '$ and 
$\la 0\ra ''$
in $N(n,k)$}\mid |T'_k|\ne 0,
|T''_k|\ne 0\rbrace\\
&\hskip1cm \ge 1-C_3n^{-C_4\varepsilon}.
\endalign
$$
\endproclaim

\demo{Proof}
Let
$$
f_\ell(s)=f^\delta_\ell(s)=\EE s^{|T^\delta_\ell|},\quad 0\le s\le 1.
$$
Then, for $0<s\le 1$
$$
\Pr\lbrace 0<|T^\delta_\ell|\le\ell\rbrace\le s^{-\ell}\EE\lbrace 
s^{|T^\delta_\ell|};|T^\delta_\ell|\ne 0\rbrace
=\frac{f_\ell(s)-f_\ell(0)}{s^\ell}.
$$
However, $\lbrace|T^\delta_\ell|\rbrace_{\ell\ge 0}$ is a
supercritical branching process, so that by Cor.\ I.11.1 in Athreya and
Ney (1972) there exists a $\lambda=\lambda(\delta)<1$ such that (in the
notation of Athreya and Ney (1972)) for all $0\le s<1$
$$
\lim_{\ell\to\infty}\frac{f_\ell(s)-f_\ell(0)}{\lambda^\ell}=Q(s)-Q(0)<\infty.
$$
Thus, if we fix $\lambda<s<1$, then there exists some $\ell_0$ such that
for $\ell\ge\ell_0$
$$
\Pr\lbrace 0<|T^\delta_\ell|\le\ell\rbrace\le 2\lbrace Q(s)-Q(0)\rbrace
\left(\frac{\lambda}{s}\right)^\ell.\tag 3.29 
$$

Next we observe that
$$
 \lim_{r\to\infty}\frac{|T^\delta_r|}{\delta^r}=W\quad \text{exists with
probability 1}, 
$$
and
$$
\Pr\lbrace W>0\rbrace=\Pr\lbrace|T^\delta_r|
\text{ is never zero}\rbrace=1-q(\delta)>0\tag 3.30
$$
(see Harris (1963) Theorems I.8.1 and I.8.3). Thus, there exists an
$\alpha=\alpha(\delta)>0$ such that
$$
\Pr\lbrace|T^\delta_r|\ge\alpha\delta^r\text{ for all }
r\ge 0\rbrace\ge\tfrac{1}{2}(1-q(\delta)).
$$
Each $\la i_1,\dots,i_\ell\ra $ in $T^\delta_\ell$ has a certain number of
descendants in $T^\delta_k$ ($k>\ell$). By the branching property these
numbers for different $\la i_1,\dots,i_\ell\ra $ are independent and have the
same distribution as $|T^\delta_{k-\ell}|$. Therefore,
$$
\align
&\Pr\lbrace |T^\delta_k|\le\alpha\delta^{k-\ell}\mid T^\delta_{[\ell]}\rbrace\\
&\hskip.5cm
  \le \Pr\lbrace\text{each individual }\la i_1,\dots,i_\ell\ra \text{ in }
T^\delta_\ell
\text{ has fewer 
 than }\alpha\delta^{k-\ell}\text{ children in }T^\delta_k 
\mid T^\delta_{[\ell]}\rbrace\\
&\hskip.5cm \le \left\lbrace 1-\tfrac{1}{2}(1-q(\delta))
  \right\rbrace^{|T^\delta_\ell|}=
\left\lbrace\frac{1+q(\delta)}{2}\right\rbrace^{|T^\delta_\ell |}.
\endalign
$$
It follows that for each $\ell<k$
$$
\align
&\Pr\lbrace 0<|T^\delta_k|<\alpha\delta^{k-\ell}\rbrace\\
&\hskip1cm \le\Pr\lbrace
0<|T^\delta_\ell|\le\ell\rbrace
+\EE\Bigl\lbrace\Pr\lbrace|T^\delta_k|<\alpha\delta^{k-\ell}\mid
T^\delta_{[\ell]}\rbrace
;|T^\delta_\ell|>\ell\Bigr\rbrace\\ 
&\hskip1cm \le  2\lbrace Q(s)-Q(0)\rbrace
\left(\frac{\lambda}{s}\right)^\ell+
\left\lbrace\frac{1+q(\delta)}{2}\right\rbrace^\ell.
\endalign
$$
If we choose $\ell\sim\varepsilon k$, then we see that for some $\beta<1$
$$
 \Pr\lbrace 0<|T^\delta_k|<\alpha\delta ^{(1-\varepsilon )k}
\rbrace\le\lbrace 2Q(s)-2Q(0)+1\rbrace\beta^{\varepsilon k}. 
$$
Since $T'_k$ and $T''_k$ are independent, each with the same distribution
as $T^\delta_k$, we conclude (use (3.30)) that
$$
\multline
\Pr\lbrace|T'_k|<\alpha\delta^{(1-\varepsilon )k}
\text{ or }|T''_k|<\alpha\delta^{(1-\varepsilon )k}\mid |T'_k|
\ne 0,|T''_k|\ne 0\rbrace\\
\le \lbrace 4Q(s)-4Q(0)+2\rbrace\lbrace1-q(\delta )
\rbrace^{-2}\beta^{\varepsilon k}.
\endmultline
$$
Finally, we observe that the conditional probability, given
$T'_{[k]},T''_{[k]}$, that there does not exist any conducting edge in
$N(n,k)$ between $T'_k$ and $T''_k$ is
$$
\left(1-\frac{\gamma (n)}{n}\right)^{|T'_k|\cdot|T''_k|}\le\exp
\left(-\frac{\gamma(n)}{n}|T'_k|\cdot|T''_k|\right).
$$
Whenever $T'_k\ne 0$, and $T''_k\ne 0$ and there is a conducting edge
between $T'_k$ and $T''_k$, then $\la 0\ra '$ 
is connected to $\la 0\ra ''$ along a
path from $\la 0\ra '$ to $T'_k$ $ (\text{in }T'_{[k]})$, then to $T''_k$ and
then to $\la 0\ra ''$ $(\text{in }T''_{[k]})$ . This path contains $(2k+1)$
edges, and since each conducting edge has resistance at most $K$ (by
(3.27)) $\rho (n,k)\le (2k+1)K$ in this situation. Thus the first equality
in (3.28) holds and 
$$
\align
&\Pr\{\text{there is no conducting path in $N(n,k)$ between
$\la 0\ra '$ and $\la 0\ra ''$}\mid |T'_k|\ne 0, |T''_k|\ne 0\rbrace\\
&\hskip1.3cm
\le \Pr\lbrace |T'_k|<\alpha\delta^{(1-\varepsilon )k}\text{ or }|T''_k|
<\alpha\delta^{(1-\varepsilon )k} \mid |T'_k|\ne 0, |T''_k|\ne 0\rbrace\\
&\hskip5cm+\exp \left(-\frac{\gamma (n)}{n} \alpha^2\delta^{2(1-\varepsilon )k}\right)\\
&\hskip1.3cm
\le\lbrace 4Q(s)-4Q(0)+2\rbrace\lbrace 1-q(\delta )\rbrace^{-2} \
\beta^{\varepsilon k}+\exp \left(-\frac{\gamma (n)}{n}\alpha^2
\delta^{2(1-\varepsilon )k}\right).
\endalign
$$
(3.28) follows when $2(1-\varepsilon )k\ge
(1-\varepsilon)(1+2\varepsilon)\log n/\log \delta$.
\qed\enddemo

For the time being we maintain the extra assumption (3.27). We use Lemma
11 to replace the random network $N(n,m_n)$ by another random network,
$M(n)$, which with high probability has a resistance at least equal to
$\rho (n,m_n)$. $M(n)$ is constructed as follows (see Figure 9). Let
$s=s_n=\lfloor\sqrt{\log n}\rfloor$.
Form $T'_{[m_n]}$ and $T''_{[m_n]}$. If $T'_s\ne 0$ and
$T''_s\ne 0$, then connect each pair of vertices
$\la i_1,\dots,i_s\ra '\in T'_s$ and $\la j_1,\dots,j_s\ra ''\in
T''_s$ which have descendants in $T'_{m_{n}}$ and $T''_{m_{n}}$,
respectively, by a resistance of size
$$
\lbrace |T'_s|+|T''_s|\rbrace\frac{K}{\log\gamma}\log n.\tag 3.31
$$
We shall write $\Cal A_{s}'$ (respectively $\Cal A_{s}''$) for the
collection of vertices $\la i_1,\dots,i_s\ra '\in T_{s}'$ (respectively
$\la j_1,\dots,j_s\ra ''\in T_{s}''$) which have descendants in
$T'_{m_{n}}$ (respectively $T''_{m_{n}}$).

\topinsert
\figure
\mletter{\la 0\ra'}{1.15}{2.3}
\mletter{T'}{2}{3.7}
\mletter{\la 0\ra''}{10.8}{2.3}
\lastletter{T''}{9.2}{3.8}
\centerline{\epsfxsize=9cm\epsffile{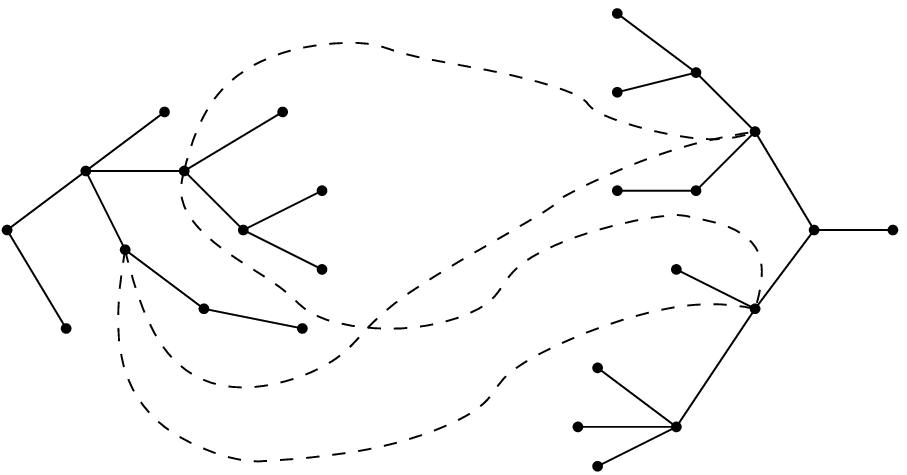}}

\capt{9}{A schematic representation of $M(n)$ with $s=2$,
$m=4$. The dashed curves represent the connections of resistance
$\lbrace|T'_s|+|T''_s|\rbrace\frac{K}{\log\gamma}\log n$.
These dashed connections have no interior points in common (when realized
in space instead of in the plane).}
\endinsert

This describes the network $M(n)$. We denote the resistance
between $\la 0\ra '$ and $\la 0\ra ''$ in $M(n)$ by 
$\Re(n,K)$ ($\Re(n,K)=\infty$ if
$T'_s$ or $T''_s$ is empty, or even if ${\Cal A}'_s$ or ${\Cal A}''_s$ is
empty).

\proclaim{Lemma 12}
If $1<\delta\le\gamma$ is such that
$$
\frac{\log\gamma}{\log\delta}<\frac{6}{5},\tag 3.32
$$
and if (3.27) holds, then
$$
\Pr\lbrace\rho(n,m_n)\le\Re(n,K)\rbrace\to 1,\tag 3.33
$$ 
and
$$
\Pr\lbrace R_n\le x\rbrace\ge \Pr\lbrace\Re(n,K)\le x\rbrace +\o(1).\tag 3.34
$$
where $\o(1)\to 0$ as $n\to\infty$.
\endproclaim

\demo{Proof}
Let $\la i_1,\dots,i_s\ra '\in T'_s$, $\la j_1,\dots,j_s\ra ''\in T_s''$
with $s=s_n=\lfloor\sqrt{\log n}\rfloor$. We apply Lemma 11 to the
sub-network of $N(n,m_n)$ consisting of the tree of descendants of
$\la i_1,\dots,i_s\ra '$ in $T'_{[m_n]}$, the tree of descendants of
$\la j_1,\dots,j_s\ra ''$ in $T''_{[m_n]}$ and the edges between the last
generations of these trees. Conditionally on $T'_{[s]},T''_{[s]}$, this
network has the same distribution as $N(n,k)$ with $k=k_n:=m_n-s_n$.
Therefore, the resistance between $\la i_1,\dots,i_s\ra '$ and
$\la j_1,\dots,j_s\ra''$ in this network has the same distribution as $\rho
(n,k_n)$. In particular, given that $\la i_1,\dots,i_s\ra '\in {\Cal
A}'_s$, $\la j_1,\dots,j_s\ra ''\in {\Cal A}''_s$, the conditional
probability that there is a path between $\la i_1,\dots,i_s\ra '$ and
$\la j_1,\dots,j_s\ra ''$ of resistance $\le (2k_n+1)K$ in the above network is
at least (for sufficiently large $n$)
$$
1-C_3n^{-C_{4}/8},\tag 3.35
$$
by virtue of (3.28) (with $\varepsilon=1/8$) and the fact that
$$
\align
k_n=m_n-s_n&\sim\frac{3}{4}\frac{\log n}{\log\gamma}\quad  (\text{see }(3.16))\\
&>\frac{5}{8}\frac{\log n}{\log\delta}\quad (\text{by }(3.32)).
\endalign
$$
Note that if the above path between $\la i_1,\dots,i_s\ra '$ and
$\la j_1,\dots,j_s\ra ''$ exists, then it is made up entirely from edges outside
$T'_{[s]}$ or $T''_{[s]}$. In fact it is built up from edges in the trees
$T'(\la i_1,\dots,i_s\ra ')$ 
and $T''(\la j_1,\dots,j_s\ra '')$ of the descendants of
$\la i_1,\dots,i_s\ra '$ and $\la j_1,\dots,j_s\ra ''$, respectively, plus an edge
between $T'_m$ and $T''_m$ which does not belong to $T'$ or to $T''$. Let
us denote by ${\Cal C}={\Cal C}_{m_n}$ the collection of conducting edges
between $T'_{m_n}$ and $T''_{m_n}$ in $N(n,m_n)$. If there is a conducting
path connecting $\la i_1,\dots,i_s\ra '$ and $\la j_1,\dots,j_s\ra ''$, then it
contains one edge from ${\Cal C}$, and this edge connects a descendant of
$\la i_1,\dots,i_s\ra '$ and a descendant of $\la j,\dots,j_s\ra ''$. Therefore, for
different pairs $\la i_1,\dots,i_s\ra '$,
$\la j_1,\dots,j_s\ra ''$ different edges from
${\Cal C}$ will be used.

Now consider the event
\itemitem {$E_n:= \{$} each pair $\la i_1,\dots,i_s\ra '\in T'$ and
$\la j_1,\dots,j_s\ra ''\in T''$ which have descendants in $T'_m$ and
$T''_m$, respectively, are connected in $N(n,m_n)$ by a conducting path in
$T'(\la i_1,\dots,i_s\ra ')\cup T''(\la j_1,\dots,j_s\ra '')\cup{\Cal C}$ of length
$(2k_n+1)$\}.

\noindent By the estimate (3.35) and the fact that ${\Cal A}'_s\subset
T'_s$, ${\Cal A}''_s\subset T''_s$ we have
$$
\Pr\lbrace E_n\mid T'_{[s]},T''_{[s]}\rbrace\ge 1-|T'_s|\cdot
|T''_s|C_3n^{-C_4/8}.
$$ 
Since (see Harris (1963) Theorem I.5.1)
$$
\EE|T'_s|=\EE|T''_s|=\delta^s=\o(n^{C_4/16}),
$$
it follows that
$$
\align
\Pr\lbrace E_n\rbrace &\ge 1-\EE\lbrace|T'_s|\cdot|T''_s|C_3n^{-C_4/8}\rbrace\\
&\ge 1-\delta^{2s}C_3n^{-C_4/8}\to 1\quad\text{as }n\to\infty.
\endalign
$$
(3.33) now follows easily from this and the monotonicity property (2.6).
Indeed for $\la i_1,\dots,i_s\ra '\in T'_s$ we may replace any edge in
$T'(\la i_1,\dots,i_s\ra ')$, with resistance $r$ say, by $|T''_s|$ parallel
edges of resistance $|T''_s|r$ without changing the resistance between any
pair of vertices in $T'_s$ and $T''_s$. Similarly we may replace any edge
in $T''_s(\la j_1,\dots,j_r\ra '')$ by $|T'_s|$ parallel edges whose resistance
is $|T'_s|$ times the resistance of the original edge. After this has been
done, we can, on the event $E_n$, connect each pair
$\la i_1,\dots,i_s\ra '\in {\Cal A}_s'$, $\la j_1,\dots,j{_s}\ra ''\in
{\Cal A}_s''$ by a path of length $(2k_n+1)$, such that the different
paths have no edges in common. Indeed each edge in $T'(\la i_1,\dots,i_s\ra ')$
has been split into $|T''_s|\ge |{\Cal A}''_s|$ parallel edges and we can
use a different one of these parallel edges to connect $\la i_1,\dots,i_s\ra '$
to different $\la j_1,\dots,j_s\ra ''$. The new paths each consist of $k_n$
edges of resistance $\le|T''_s|K$, an edge of ${\Cal C}$ of resistance
$\le K$ and $k_n$ edges of resistance $\le|T'_s|K$, all of these edges
being in series. The resistance of such a path is therefore at most
$$
\lbrace k_n(|T'_s|+|T''_s|)+1\rbrace K\le\lbrace|T'_s|+|T''_s|
\rbrace\frac{K}{\log\gamma}\log n.
$$ 
The new paths between all the pairs $\la i_1,\dots,i_s\ra '\in{\Cal
A}_s''$ and $\la j_1,\dots,j_s\ra ''\in{\Cal A}_s''$ are edge-disjoint,
but they still have vertices in common in $T'_{[m_n]}$ and in
$T''_{[m_n]}$. However, by (2.6) these contacts between different paths can
only reduce the resistance between $\la 0\ra '$ and $\la 0\ra ''$. Therefore, on
$E_n$, the resistance between $\la 0\ra '$ and $\la 0\ra ''$ in $N(n,m_n)$ is at most
the resistance between $\la 0\ra '$ and $\la 0\ra ''$ in $M(n)$. This proves (3.33).
(3.34) follows from (3.33) and (3.26).
\qed\enddemo

Apart from removing some truncations later on, the only estimate left is
one which shows that $\Re(n,K)$ is essentially equal to the sum of the
resistance of $T'_{[s]}$ and $T''_{[s]}$. This is done in the next lemma
by showing that all vertices in ${\Cal A}_s'$ and ${\Cal A}_s''$ have
almost the same potential.

\proclaim {Lemma 13}
Assume that for some $0<\varepsilon<K<\infty$
$$
F(\varepsilon-)=0,\quad F(K)=1.\tag 3.36
$$
Then
$$
\Pr\lbrace\Re(n,K)\le x\rbrace\to\Pr\lbrace R'(\delta)+R''(\delta)\le
x\rbrace\tag 3.37
$$
at each continuity point of the right hand side.
\endproclaim

\demo{Proof}
We shall only need
$$
 \liminf_{n\to\infty}\Pr\lbrace\Re(n,K)\le x\rbrace\ge\Pr\lbrace
R'(\delta)+R''(\delta )\le x\rbrace.\tag 3.38 
$$
We therefore only prove (3.38) and leave the (easy) other half of (3.37)
to the interested reader. First we note that $\Re(n,K)=\infty$ if ${\Cal
A}'_s=\emptyset$ or ${\Cal A}''_s=\emptyset$, or equivalently, if
$T'_{m_n}=\emptyset$ or $T''_{m_n}=\emptyset$. But as $n\to\infty$ also
$m_n\to\infty$ and
$$
\align
\Pr\lbrace T'_{m_n}=\emptyset\rbrace=\Pr\lbrace T''_{m_n}=\emptyset 
\rbrace&\to\Pr\lbrace T^\delta\text{ is finite}\rbrace=
\Pr\lbrace Z\text{ dies out}\rbrace\\
&=q(\delta)=\Pr\lbrace R'(\delta)=\infty\rbrace=\Pr\lbrace 
R''(\delta)=\infty\rbrace,
\endalign
$$
where $Z$ is the branching process corresponding to $T^\delta$ (apply Lemma
2 to $T^\delta$). We therefore should show (at continuity points $x$ of
the right hand side of (3.37))
$$
\align
\liminf_{n\to\infty}\Pr &\lbrace\Re(n,K)\le x\mid {\Cal
A}'_s\ne\emptyset,{\Cal A}_s''\ne\emptyset\rbrace\tag 3.39\\ 
&\ge\Pr\lbrace
R'(\delta)+R''(\delta)\le x\mid T'\text{ and }T''\text{ are
infinite}\rbrace.
\endalign
$$
Now assume that ${\Cal A}_s'\ne\emptyset$, ${\Cal A}_s''\ne\emptyset$ and
let $\lbrace X_\nu\rbrace_{\nu\ge 0}$ be a Markov chain on the vertices of
$T'_{[s_n]}\cup T_{[s_n]}''$ with transition probability matrix
$$
P(\la y\ra ,\la z\ra )=
\left\lbrace\sum_{y}\frac{1}{R(e)}\right\rbrace^{-1}
\frac{1}{R(\la y\ra ,\la z\ra)}\tag
3.40
$$ 
whenever $\la y\ra $ and $\la z\ra $ are neighbors in $M(n)$ (by this we mean that
$\la y\ra $ and $\la z\ra $ both belong to $T_{[s_n]}'$ or both to $T_{[s_n]}''$, and
are neighbors in $T_{[s_n]}$, respectively $T_{[s_n]}''$, or one belongs
to ${\Cal A}_{s_n}'$ and the other to ${\Cal A}_{s_n}''$). In (3.40)
$R(\la y\ra ,\la z\ra )$ denotes the resistance of the edge between $\la y\ra $ and $\la z\ra $
(if $\la y\ra $ and $\la z\ra $ are neighbors), and $\sum_y$ runs over all edges of
$M(n)$ with one endpoint at $y$. By (2.3) and (2.4) $\Re(n,K)$, the
resistance between $\la 0\ra '$ and $\la 0\ra ''$ in $M(n)$, equals
$$
\biggl\lbrace\sum_{\la i_1\ra \in T'} \frac{1}{R(\la 0\ra ',\la i_1\ra ')}
\Pr\lbrace X.\text{ reaches }\la 0\ra ''\text{ before }
\la 0\ra '\mid X_0=\la i_1\ra '\rbrace\biggr\rbrace^{-1}.
\tag "\rlap{(3.41)}"
$$
Furthermore, it is probabilistically evident that
$$
\align
&\Pr\lbrace X.\text{ reaches }\la 0\ra ''\text{ before }\la 0\ra '
\mid X_0=\la i_1\ra '\rbrace
\tag3.42\\ 
&\hskip.5cm=\sum_{\la x\ra \in{\Cal A}_s'}\Pr\lbrace X.\text{ reaches }
{\Cal A}_s'\text{ before }\la 0\ra '\text{ and does so 
 first at }\la x\ra '\mid X_0=\la i_1\ra '\rbrace\\
&\hskip3.8cm\times\Pr\lbrace X.\text{ reaches }\la 0\ra ''
\text{ before }\la 0\ra '\mid X_0=\la x\ra '\rbrace.
\endalign
$$
Assume that we can prove the existence of a (random) sequence of numbers
$P'_n$ such that
$$
\sup_{\la x\ra '\in{\Cal A}{_s}'}| \Pr\lbrace X.\text{ reaches }
\la 0\ra ''\text{ before }\la 0\ra '\mid X_0=\la x\ra '\rbrace
-P_n'|\to 0\quad (n\to\infty)\tag 3.43
$$
in probability on the event
$$
F_n:=\lbrace{\Cal A}_{s_n}'\ne\emptyset,{\Cal A}_{s_n}''
\ne\emptyset\rbrace.\tag 3.44 
$$
(Of course this means that the probability of the subset of $F_n$ on which
(3.43) fails tends to 0.) Then (3.41)--(3.43) yield
$$
\align
&\Re^{-1}(n,K)=(P'_n+\o_n(1))\tag 3.45\\
&\hskip2cm\times\sum_{\la i_1\ra '\in T'}
\frac{1}{R(\la 0\ra ',\la i_1\ra ')}\Pr\lbrace X.
\text{ reaches }{\Cal A}{_s}'\text{ before }\la 0\ra '\mid X_0=\la i_1\ra '\rbrace
\endalign
$$
where $\o_n(1)\to 0$ in probability on $F_n$ as $n\to\infty$. Again by
(2.3) and (2.4) the sum in the right hand side of (3.45) equals the
reciprocal of the resistance between $\la 0\ra '$ and ${\Cal A}_s'$ in
$T_{[s]}'$. For the time being we denote the latter resistance by $\Re
_{n}'$. With this notation (3.45) can be written as \footnote{Note that
$\Re_n'$ and $\Re(n,K)$ cannot be zero under (3.36).}
$$
P'_n-\frac{\Re'_n}{\Re(n,K)}\to 0\quad\text{in probability on }F_n.\tag 3.46
$$
If we can prove (3.43), then by interchanging the roles of $T'$ and $T''$
we can also prove the existence of a $P_n''$ such that
$$
\sup_{\la y\ra ''\in{\Cal A}_s''}|\Pr\lbrace X.\text{ reaches }\la 0\ra '
\text{ before }\la 0\ra ''\mid
X_0=\la y\ra ''\rbrace-P''_n|\to 0\tag 3.47
$$  
in probability on $F_n$, and
$$
P_n''-\frac{\Re''_n}{\Re(n,K)}\to 0\quad\text{in probability on }F_n,\tag 3.48
$$
where $\Re_n''$ is the resistance between $\la 0\ra ''$ and ${\Cal A}_s''$
in $T_{[s]}''$.
Furthermore, (3.43) and (3.47) together imply on $F_n$, uniformly in $\la  x
\ra ' \in {\Cal A}'_s$,
$$\align
&\Pr\{X_.\text{ reaches }\la 0\ra''\text{ before }\la 0\ra'\mid 
  X_0 = \la  x \ra '\}\tag 3.49\\
&\hskip.5cm= P'_n + \o_n(1)\\
&\hskip.5cm= \sum_{\la  y\ra '' \in {\Cal A}''_s} \Pr\{X_. 
\text{ reaches } {\Cal A}''_s \text{ before }
\la  0 \ra ' \text{ and does so first at } 
  \la  y \ra '' \mid X_0 = \la  x \ra '\}\\
&\hskip1.5cm\times \Pr 
\{X. 
\text{ reaches } \la  0\ra '' \text{ before } 
   \la  0\ra ' \mid X_0 = \la  y \ra ''\} \\
&\hskip.5cm= \Pr\{ X_. \text{ reaches } {\Cal A}_s'' \text{ before } \la  0\ra ' 
\mid X_0 = \la  x'\ra  \} (1-P_n'' + \o_n(1)).
\endalign
$$

We shall next prove that
$$
\min_{\la  x \ra ' \in {\Cal A}'_{s}} \Pr \{X_. \text{ reaches } 
{\Cal A}''_s \text{ before } \la  0 \ra ' \mid X_0 = \la  x \ra '\} \to 1 
\quad\text{in probability on } F_n . \tag "\rlap{(3.50)}"
$$
(Once one has (3.43) and (3.47) it is not hard to obtain (3.50) by
equating the current flowing into ${\Cal A}'_s$ and the current flowing
out of ${\Cal A}'_s$.  However, we need (3.50) to prove (3.43) so we must
prove it directly.) The construction of $M(n)$ is such that the resistance
of any edge between any pair of vertices $\la  u\ra '\in {\Cal A}'_s$
and $\la v\ra''\in\Cal A_s''$ has the same value, namely the value
in (3.31). But there are $|\Cal A_s''|$ edges between any $\la u\ra'\in\Cal A_s'$ and the
set ${\Cal A}''_s$, and the only other edge incident to $\la  u \ra '$
is an edge from $\la  u \ra '$ to $T'_{s-1}$ with resistance between
$\varepsilon$ and $K$ (by (3.36)).  Therefore for any $\la  u \ra ' \in
{\Cal A}'_s$.
$$
P\{ \la u\ra ',\la v\ra ''\}\text{ has the same value for all }\la v\ra ''
\in{\Cal A}''_s,\tag 3.51
$$
and the value in (3.51) satisfies
$$
\sum_{\la v\ra ''\in{\Cal A}''_s} P\{ \la u\ra ', \la v\ra ''\}
\ge\frac{|{\Cal A}''_s|\log \gamma}{K\{|T'_s|+|T''_s|\}\log n}
\bigg\{\frac{|{\Cal A}''_s|\log\gamma}{K\{|T'_s|+|T''_s|\}\log n}
+\frac{1}{\varepsilon}\bigg\}^{-1}.\tag 3.52
$$
A decomposition with respect to the last visit of $X_.$ to ${\Cal A}'_s$
before it hits $\la 0\ra '$ or ${\Cal A}''_s$ yields.
$$
\align
&\Pr\{ X.\text{ reaches } \la 0\ra '\text{ before }{\Cal A}''_s
\mid X_0=\la x\ra '\}\tag 3.53\\
&\hskip.5cm= \sum_{\la u\ra '\in{\Cal A}'_s} \EE\{\text{number of visits to } 
\la u\ra '\text{ before
 reaching }\la 0\ra '\text{ or }{\Cal A}''_s\mid
X_0=\la x\ra '\}\\
&\hskip2cm  
\times\Pr\{ X_.\text{ reaches } 
\la 0\ra ' \text{ without returning to }{\Cal A}'_s\mid X_0=\la u\ra '\}.
\endalign
$$
Similarly,
$$\align
1 &= \Pr\{ X_.\text{ reaches }\la 0\ra '\text{ or }{\Cal A}''_s
\text{ sometime}\mid X_0=\la x\ra '\}\tag 3.54\\
  &= \sum_{\la u\ra '\in{\Cal A}'_s} \EE\{
\text{number of visits to } \la u\ra '\text{ before
 reaching } \la 0\ra '\text{ or }{\Cal A}''_s\mid X_0=\la x\ra '\}\\
&\hskip1cm\times \Pr\{ X_.
\text{ reaches } \la 0\ra '\text{ or }{\Cal A}''_s
\text{  without returning to }{\Cal A}'_s\mid X_0=\la u\ra '\}.
\endalign
$$
Dividing (3.53) by (3.54) we see that
$$
\align
&\Pr\{ X_.\text{ reaches }\la 0\ra '\text{ before }{\Cal A}''_s
\mid X_0=\la x\ra '\}\tag 3.55\\
&\le\max_{\la u\ra '\in{\Cal A}'_s}\frac{\Pr\{ X_.\text{ reaches }\la 0\ra '\text
{ before }{\Cal A}'_s\mid X_0=\la u\ra '\}}{\Pr\{ X_.\text{ reaches }\la 0\ra '\text{ or }
{\Cal A}''_s\text{ before returning to }{\Cal A}'_s\mid X_0=\la u\ra '\}}.
\endalign
$$

To estimate (3.55) note first that
$$
{\Cal A}''_s\subset T''_s;\text{ hence }|{\Cal A}''_s|\le|T''_s|.\tag 3.56
$$
Also, there exist random variables $W'$, $W''$ such that
$$
{|T_s'|\over\delta^s}\to W',\ {|T_s''|\over\delta^s}\to W''\quad \text{w.p.1},
$$
and
$$
\align
W'&>0 \text{ a.e., on the set }\{ T'_p\neq\emptyset\text{ for all }p\}\\
\text{and }W''&>0 \text{ a.e., on the set }\{ T''_p\neq\emptyset\text{ for all }p\}
\endalign
$$
(see  Harris (1963) Theorem I.8.1 and Remark I.8.1).
Thus $|T'_s|$ and $|T''_s|$ are both of order $\delta^s$ on most of the set
$$
\{ T'_s\ne\emptyset,T''_s\ne\emptyset\}.
$$
Moreover, by definition of ${\Cal A}'_s$,
$$\align
\Pr\{\la u\ra '\in{\Cal A}'_s\mid T'_{[s]}, \la u\ra '\in
T'_s\} &= \Pr\{ \la u\ra '\text{ has descendants in
}T'_{[m]}\mid \la u\ra '\in T'_s\}\\ 
&=\Pr\{ T'_{m-s}\ne\emptyset\}\\
&\to
\Pr\{ T'_p\ne\emptyset\text{ for all }p\} =1-q(\delta)>0
\endalign
$$
(compare (2.8) and (2.9)). If $T'_s=\{ \la u_1 \ra ',\dots,\la u_t\ra '\}$, then by
the branching property, the events $\la u_i\ra '\in{\Cal
A}'_s$, $i=1,\dots,t$, are conditionally independent, given
$\la u_i\ra '\in T'_s$, $1\le i\le t$. It follows from these
observations that
$$
I[T'_s\ne\emptyset ]\bigg\{ {|{\Cal A}_s|\over |T'_s|}-(1-q(\delta))
\bigg\}\to 0\quad\text{in probability}.
$$
The same relation holds when ${\Cal A}'_s$ and $T'_s$ are replaced by
${\Cal A}''_s$ and $T''_s$. Since ${\Cal A}''_s\ne\emptyset$ implies
$T''_s\ne\emptyset$, and similarly for ${\Cal A}'_s$ and $T'_s$ we find
also that (cf.\ (3.56)).
$$
{|{\Cal A}'_s|\over |T'_s|}\text{ and }{|{\Cal A}''_s|\over |T''_s|}\to 
1-q(\delta)\quad\text{in probability on the }
\text{set }F_n\text{ of (3.44)}
$$
and
$$
{|{\Cal A}''_s|\over |T'_s|+|T''_s|}\to {(1-q(\delta))W''\over W'+W''}
\quad \text{in probability on 
the set }F_n.\tag 3.57
$$

We return to (3.55). First we estimate its denominator.
$$
\align
&\min_{\la u\ra '\in{\Cal A}'_s}\quad \Pr\{ X_.\text{ reaches }\la 0\ra '
\text{ or }
{\Cal A}''_s\text{ before returning to }
{\Cal A}'_s\mid X_0=\la u\ra '\}\tag 3.58\\
&\hskip2cm\ge\sum_{\la v\ra ''\in{\Cal A}''_s} 
P\{\la u\ra ',\la v\ra ''\}.
\endalign
$$
By virtue of (3.52) and (3.57), for every $\eta>0$ and all large $n$ there
exists a $\kappa(\eta)>0$ such that the subset of $F_n$ on which the right
hand side of (3.58) is less than
$$
\frac{\varepsilon}{\kappa(\eta)\log n}
$$
has probability $\le\eta$. On the other hand, for all
$\la u\ra '\in{\Cal A}'_s$ the numerator in the right hand side of
(3.55) is bounded above by
$$
\sup_{w\in T'_{s-1}}\pi(\la w\ra ', T'_{[m]},R,s)=\Pi(T'_{[m]},R,s)
$$
(by (2.49)). Therefore, by Lemma 6, the numerator in the right hand 
side of (3.55) is at most (for large $n$)
$$
\bigg(\frac{2L}{2L+\varepsilon}\bigg)^{C_1s_n}
$$
on the set $F_n$ minus a subset of probability at most $\exp(-C_2s_n)$.
Consequently, for large $n$, on the set $F_n$, minus a set of probability
at most $2\eta$, (3.55) is at most
$$
\bigg(\frac{2L}{2L+\varepsilon}\bigg)^{C_1s_n}\quad
\frac{\kappa (\eta)\log n}{\varepsilon}\to 0
$$
(as $n\to\infty$, since $s_n=\sqrt{\log n})$. Since $\eta$ is arbitrary we
finally proved that the right hand side of (3.55) tends to $0$ in
probability on $F_n$. This implies (3.50).

From here on it is easy to complete the proof of the lemma (still under
the assumptions (3.43) and (3.47)). Firstly, (3.49) and (3.50) together
imply
$$
P'_n+P''_n\to 1\quad \text{in probability on }F_n.\tag 3.59
$$
In turn (3.59), (3.46), and (3.48) together show
$$
\frac{\Re'_n+\Re''_n}{\Re(n,K)}\to 1\quad\text{in probability on }F_n.\tag 3.60
$$ 
It follows from (2.3), (2.4) and the definition ${\Cal A}'_s$ that
$R(T'_{[s_n]})\le\Re'\le R(T'_{[m_n]})$. A similar inequality holds for
$\Re''$ so that
$$
\Re'_n\to R(T'),\  \Re''_n\to R(T'')\quad\text{w.p.1.}
$$
Since $R(T')$ and $R(T'')$ both have the distribution of $R(T^\delta)$,
that is of $R'(\delta)$ and $R''(\delta)$, (3.60) implies (3.39) and the
lemma. The proof has therefore been reduced to proving (3.43) and (3.47).

(3.43) can now be proved quickly from (3.52) and (3.51). Indeed when $X_.$
starts at $\la x\ra '\in{\Cal A}'_s$ it cannot reach $\la 0\ra ''$
without passing through ${\Cal A}''_s$. Let $\tau$ be the first time $X_.$
visits ${\Cal A}''_s$. Then a decomposition with respect to the values of
$\tau-1$, $X_{\tau-1}$ and $X_\tau$ gives
$$\align
&\Pr\{ X_.\text{ reaches }\la 0\ra ''\text{ before }\la 0\ra ' 
  \mid X_0=\la x\ra '\}\tag 3.61\\
&\hskip1cm = \sum^\infty_{n=1}\sum_{\la u\ra '\in{\Cal A}'_s}
\sum_{\la v\ra ''\in{\Cal A}''_s} \Pr\{ X_n=\la u\ra ',X_j\ne \la 0\ra ',
X_j\notin{\Cal A}''_s, 0\le j\le n\}\\
&\hskip2cm\times P\{ \la u\ra ',\la v\ra ''\} \Pr\{ X_.\text{ reaches }
\la 0\ra ''\text{ before }\la 0\ra '\mid X_0=\la v\ra ''\}.
\endalign
$$
By virtue of (3.51)
$$
\align
&\sum_{\la v\ra ''\in{\Cal A}''_s}P\{ \la u\ra ',\la v\ra ''\} \Pr\{ X_.
\text{ reaches }\la 0\ra ''\text{ before }
\la 0\ra '\mid X_0=\la v\ra'' \}\\
&\hskip1cm  =\frac{1}{|{\Cal A}''_s|}\sum_{\la v\ra ''\in
{\Cal A}''_s}P\{\la u\ra ',\la v\ra ''\}\\
&\hskip2cm\times
\sum_{\la w\ra ''\in{\Cal A}''_s}\Pr\{ X_.\text{ reaches }\la 0\ra ''
\text{ before }\la 0\ra '\mid X_0=\la w\ra ''\}.
\endalign
$$
Thus, (3.61) equals
$$
\align
&\Pr\{ X_.\text{ reaches }{\Cal A}''_s\text{ before }\la 0\ra '
\mid X_0=\la x\ra '\}\\
&\hskip2cm\times\frac{1}{|{\Cal A}''_s|}\sum_{\la w\ra ''\in{\Cal A}''_s}
\Pr\{ X_.\text{ reaches }\la 0\ra ''\text{ before }\la 0\ra '
\mid X_0=\la w\ra ''\},
\endalign
$$
which, together with (3.50) implies (3.43) with
$$
P'_n=\frac{1}{|{\Cal A}''_s|}\sum_{\la w\ra ''\in{\Cal A}''_s}
\Pr\{ X_.\text{ reaches }\la 0\ra ''\text{ before }
\la 0\ra '\mid X_0=\la w\ra ''\}.
$$
This proves (3.43), and as observed before, (3.47) follows by
interchanging the roles $T'$ and $T''$ in the proof of (3.43).
\qed\enddemo

 As a result of Lemmas 12 and 13 we have
$$
 \liminf_{n\to\infty} \Pr\{ R_n\le x\}\ge \Pr\{
R'(\delta)+R''(\delta)\le x\}\tag 3.62 
$$
at each continuity point of the right hand side, whenever (3.32) and
(3.36) hold. Now set
$$
 R(e,\varepsilon,K)=\cases R^\varepsilon (e)=R(e)+\varepsilon
&\text{if $R^\varepsilon\le K$}\\ 
\infty&\text{if
$R^\varepsilon>K$}.\endcases\tag 3.63 
$$
When $R(e)$ is replaced by $R(e,\varepsilon ,K)$ we shall write
$R_n(\varepsilon ,K)$ (respectively $R(\delta ;\varepsilon ,K)$ or
$R(T^\delta ;\varepsilon ,K)$) for the resistance between $0$ and $\infty$
in $K_{n+2}$ (respectively between $\la 0\ra $ and $\infty$ in $T^\delta$).
(3.62) applies when $R(e)$ is replaced by $R(e,\varepsilon,K)$ (provided
(3.32) holds). This replacement only increases resistances. Consequently
$$
\align
\liminf_{n\to\infty}\Pr\{ R_n\le x\}&\ge\liminf_{n\to\infty}
\Pr
\{R_n(\varepsilon, K)\le x\}\\
&\ge\Pr\{ R'(\delta;\varepsilon,K)+R''(\delta;\varepsilon,K)\le x\} 
\endalign
$$
for each $0<\varepsilon<K<\infty$ and $\delta$ satisfying (3.32). Here, of
course, $R'(\delta;\varepsilon,K)$ and $R''(\delta;\varepsilon,K)$ are
independent copies of $R(\delta;\varepsilon,k)$. It follows that
$$
\liminf_{n\to\infty} \Pr\{ R_n\le x\}\ge
\lim_{\varepsilon\downarrow 0}\lim_{\delta\uparrow\gamma}
\lim_{K\to\infty} \Pr\{ R'(\delta;\varepsilon,K)
+R''(\delta;\varepsilon,K)\le x\}.
$$
Theorem 3 will therefore be a consequence of Lemma 8 and the next lemma in
which we remove the ``truncations".

\proclaim{Lemma 14} As $K\to\infty,$ $\delta\uparrow\gamma$ and
$\varepsilon\downarrow 0$ (in this order)
$R(\delta;\varepsilon,K)$ converges in distribution to
$R(\gamma)=R(T^\gamma)$.
\endproclaim

\demo{Proof} For any realization $t,r(\cdot)$ of $T^\delta$ and $R(\cdot)$
we consider a Markov chain $\{X_\nu\}=\{X_\nu(\varepsilon,K,t,r)\}$ on $t$
with transition probability matrix
$$\align
P(\la y\ra ,\la z\ra )&=P(\la y\ra ,\la z\ra ;\varepsilon,K,t,r)
\tag "\rlap{(3.64)}"\\
&=\biggl\{\sum _y{1\over r(e;\varepsilon,K)}\biggr\}^{-1}
{1\over r(y,z;\varepsilon,K)}, 
\quad y,z\text{ adjacent on }t,
\endalign
$$
where again $\sum_y$ runs over all edges $e$ incident to $y$,
$r(e;\varepsilon,K)$ is defined by (3.63) with $R$ replaced by $r$ and
$r(y,z;\varepsilon,K)$ is $r(e;\varepsilon,K)$ for $e$ the edge between
$y$ and $z$. As before we take $P(\la y\ra ,\la z\ra )=0$ if $y$ and $z$ are not
adjacent in $t$. By (2.3) and (2.4)
$$
\multline
\{R (\delta;\varepsilon,K)\}^{-1}\\
=\sum_{\la i\ra \in T^\delta_1}
\frac{1} {R(e(i);\varepsilon,K)}\Pr\{ X_\nu(\varepsilon,K,T^\delta,R)
\text{ reaches }\infty\text{ before }\la 0\ra \mid X_0=\la i\ra \}.
\endmultline
$$
First we show that we can take $K=\infty$, i.e., we prove
$$
\align
&\Pr\{ X_\nu(\varepsilon,K,T^\delta,R)\text{ reaches }\infty\text{ before }
\la 0\ra \mid X_0=\la i\ra \}\tag 3.65\\
&\hskip2cm\to \Pr\{ X_\nu(\varepsilon,\infty,T^\delta,R)\text{ reaches }\infty
\text{ before }\la 0\ra \mid X_0=\la i\ra \}
\endalign
$$
in probability as $K\to\infty$. Here
$X_\nu(\varepsilon,\infty,T^\delta,R)$ is defined by (3.64) and (3.63)
with $K=\infty$, i.e., the resistance of $e$ is taken as
$R^\varepsilon(e)$. Note that when we refer to ``convergence in
probability'' in (3.65) we view both sides as random variables, namely as
functions of $T^\delta$ and $R(\cdot)$. Of course both sides are zero when
$T^\delta$ is finite so that we can restrict ourselves to that part of the
probability space where $T^\delta$ is infinite. On this event we define
$$
\Cal B_s=\Cal B_s^\delta=\{\la x\ra \in T^\delta_s:T^\delta(x),
\text{ the tree of descendants
 of }\la x\ra,\text{ is infinite}\}.
$$
Clearly for any realization $t$, $r$ of $T$, $R$ with $t$ infinite
$$
\multline
\Pr\{ X_\nu(\varepsilon,K,t,r)\text{ reaches }\infty
\text{ before }\la 0\ra \mid X_0=\la i\ra \}\\
\le\Pr\{ X_\nu(\varepsilon,K,t,r)\text{ reaches }\Cal B_s
\text{ before }\la 0\ra \mid X_0=\la i\ra \}.
\endmultline
$$
Also, for $s\ge 1$, $X_\nu$ is contained in a finite graph until it reaches
$\Cal B_s$, so that
$$
\multline
\lim_{K\to\infty}\Pr\{ X_\nu(\varepsilon,K,t,r)\text{ reaches }
\Cal B_s\text{ before }\la 0\ra \mid X_0=\la i\ra \}\\
=\Pr\{ X_\nu(\varepsilon,\infty,t,r)\text{ reaches }
\Cal B_s\text{ before }\la 0\ra \mid X_0=\la i\ra \}.
\endmultline
$$
Therefore
$$
\align
&\limsup_{K\to\infty}\Pr\{ X_\nu(\varepsilon,K,T^\delta,r)
\text{ reaches }\infty\text{ before }\la 0\ra \mid X_0=\la i\ra\}
\tag "\rlap{(3.66)}"\\
&\hskip1cm\le\lim_{s\to\infty}\Pr\{ X_\nu(\varepsilon,\infty,T^\delta,r)
\text{ reaches }\Cal B_s\text{ before } \la 0\ra \mid X_0=\la i\ra \}\\
&\hskip1cm=\Pr\{ X_\nu(\varepsilon,\infty,T^\delta,R)\text{ reaches }\infty
\text{ before }\la 0\ra \mid X_0=\la i\ra \}.
\endalign
$$
On the other hand
$$
\align
&\Pr\{ X_\nu(\varepsilon,K,T^\delta,R)\text{ reaches }\infty
\text{ before }\la 0\ra \mid X_0=\la i\ra \}\tag 3.67\\
&\hskip1cm=\sum_{\la x\ra \in\Cal B_s}\Pr\{ X_\nu(\varepsilon,K,T^\delta,R)
\text
{ reaches }\Cal B_s\text{ before }\la 0\ra \text{ and }\\
&\hskip4cm\text{ reaches }\Cal B_s\text{ first at }\la x\ra\mid X_0=\la i\ra\}\\
&\hskip3cm\times
\Pr\{ X_\nu(\varepsilon,K,T^\delta,R)\text{ reaches }\infty\text
{ before }\la 0\ra \mid X_0=x\}.
\endalign
$$
By Lemma 3 the last factor under the sum in the right hand side of (3.67)
is at least
$$
\frac{\rho(x;\varepsilon,K)}{\rho(x;\varepsilon,K)+R(T^\delta(x);
\varepsilon,K)},
\tag 3.68
$$
where for $\la x\ra =\la i_1,\dots,i_s\ra \in T^\delta_s$
$$
\rho(x;\varepsilon,K)=\sum_{j=1}^s R(e(i_1,\dots,i_j);\varepsilon,K)
\ge s\varepsilon.
$$
Consequently (3.68) is at least
$$
1-\min\left\{\frac{R(T^\delta(x);\varepsilon,K)}{s\varepsilon},1\right\},
$$
and for each $\eta\in (0,1)$
$$
\align
&\Pr\{ X_\nu(\varepsilon,K,T^\delta,R)\text{ reaches }\infty
\text{ before }\la 0\ra \mid X_0=\la i\ra \}\\
&\hskip1cm\ge (1-\eta)\Pr\{ X_\nu(\varepsilon,K,T^\delta,R)\text{ reaches }\Cal B_s
\text{ before }\la 0\ra \mid X_0=\la i\ra \}\\
&\hskip4cm-\Pr\{ R(T^\delta(Y_s);\varepsilon,K)\ge\eta\varepsilon s\mid X_0=\la i\ra \},
\endalign
$$
where $Y_s$ denotes the position where $X_.$ hits $\Cal B_s$ first,
provided $X_0$ does hit $\Cal B_s$; if $X_.$ does not hit $\Cal B_s$ then
$R(T^\delta(Y_s);\varepsilon,K)$ is taken to be zero. But until $X_.$ hits
$\Cal B_s$ it cannot enter $T^\delta(Y_s)$ and so knowledge of $Y_s$
contains no information on $T^\delta(Y_s)$, nor on its resistances, except
that $T^\delta(Y_s)$ must be infinite (that is the meaning of
$Y_s\in\Cal B_s)$. Therefore the conditional distribution of
$R(T^\delta(Y_s);\varepsilon,K)$ given $T^\delta_{[s]}$,
$\Cal B_s$ and $Y_s$,
is simply the conditional distribution of $R(T^\delta;\varepsilon,K)$
given that $T^\delta$ is infinite. Thus, to prove (3.65) we merely have to
prove that for each $\sigma>0$ there exists a $K_\sigma<\infty$ and
$y_\sigma<\infty$ such that
$$
\Pr\{ R(T^\delta ;\varepsilon,K)\ge y\mid T^\delta\text{ is infinite}\}
\le 2\sigma
\quad\text{for all } K\ge K_\sigma\text{ and }y\ge y_\sigma.\tag 3.69
$$
To prove (3.69) first choose $\widehat K$ such that
$$
\delta \Pr\{ R(e;\varepsilon,\widehat K)<\infty\}=\delta F(\widehat K-
\varepsilon)>1.
$$
Now consider a $K>\widehat K.\quad R(T^\delta;\varepsilon,K)$ is the
resistance between $\la 0\ra $ and $\infty$ in the tree obtained by removing
from $T^\delta$ each edge $e$ with $R(e)>K-\varepsilon$ and by replacing
$R(e)$ by $R^\varepsilon(e)$ on the other edges. If with an edge from
$\la i_0,\dots,i_n\ra $ to $\la i_0,\dots,i_{n+1}\ra $ which gets removed 
--- because
its resistance exceeds $(K-\varepsilon)$ --- we also remove all of
$T^\delta(i_0,\dots,i_{n+1})$, then the resulting tree is just the family
tree $\widetilde T$ of a branching process whose offspring distribution
puts mass
$$
\widetilde{p}_m:=\sum_{n\ge m}\Pr\{|T_1^\delta|=n\}
\binom nm F^m(K-\varepsilon)(1-F(K-\varepsilon))^{n-m}
$$
on $m$ (compare with $\widetilde T$ and $\widetilde{p}_m$ in Lemma 2). The
mean of this offspring distribution is
$$
\EE|T^\delta_1|F(K-\varepsilon)=\delta F(K-\varepsilon)\ge\delta F
(\widehat K-\varepsilon)>1.
$$
Consequently, by Lemma 2
$$
\Pr\{ R(T^\delta;\varepsilon,K)=\infty\}=\Pr\{\widetilde T
\text{ is finite}\}=\widetilde{q}(K),\tag 3.70 
$$
where $\widetilde{q}=\widetilde{q}(K)$ is the unique root in $[0,1)$ of 
$$
x=\EE x^{|\widetilde T_1|}=\sum^\infty_{m=0}\widetilde{p}_m x^m.
$$
Now, since $\widetilde T$ is a subgraph of $T^\delta$,
$$\align
&\Pr\{ R(T^\delta;\varepsilon,K)\ge y\mid T^\delta\text{ is infinite}\}
\tag 3.71\\
&\hskip.5cm
= {1\over 1-q(\delta)}\Pr\{ R^\varepsilon(\widetilde T)\ge y\text{ and }
T^\delta\text{ is infinite}\}\\
&\hskip.5cm 
= {1\over 1-q(\delta)}\Bigl[\Pr\{\widetilde T\text{ is finite but }T^\delta
\text{ is infinite}\}
+ \Pr\{ R^\varepsilon(\widetilde T)\ge y\text{ and }\widetilde T
\text{ is infinite}\}\Bigr]\\
&\hskip.5cm
= {\widetilde{q}(K)-q(\delta)\over 1-q(\delta)}+{1-
\widetilde{q}(K)\over 1-q(\delta)}\Pr\{ R^\varepsilon(\widetilde T)\ge 
y\mid \widetilde T\text{ is infinite}\rbrace.
\endalign
$$
Since
$$
\lim_{K\to\infty} \EE  x^{\v \widetilde T_{1}\v } = 
\EE  x^{\v T^{\delta}_{1}\v}
\quad\text{uniformly on } [0,1],
$$
one easily sees that
$$
\lim_{K\to\infty}\widetilde{q}(K)=q(\delta).
$$
We can therefore first choose $K_\sigma>\widehat K$ large to make
$$
\frac{\widetilde{q}(K_\sigma)-q(\delta)}{1-q(\delta)}\le\sigma,
$$
and then choose $y$ so large that the second term in the last member of
(3.71) with $K=K_\sigma$ is at most $\sigma$ (by (3.70)). Since
$R(T^\delta;\varepsilon,K)\le R(T^\delta;\varepsilon,K_\sigma)$ for $K\ge
K_\sigma$ (by (2.6)) this proves (3.69), and consequently also (3.65).

(3.65) allows us to take $K=\infty$. Next we should prove 
$$
\align
&\Pr\{ X_\nu(\varepsilon,\infty,T^\delta,R)\text{ reaches }\infty
\text{ before }\la 0\ra \mid X_0=\la i\ra \}\tag 3.72\\
&\hskip2cm\to \Pr\{ X_\nu(\varepsilon,\infty,T^\gamma,R)\text{ reaches }\infty
\text{ before }\la 0\ra \mid X_0=\la i\ra \}
\endalign 
$$
in probability as $\delta\uparrow\gamma$. This will allow us to prove that
we may take $\delta=\gamma$. We do not give the details for (3.72). It is
very similar to the proof of (3.65) if we take account of the fact that we
realize $T^\delta$ for all $\delta<\gamma$ and $T^\gamma$ on the same
probability space such that $T^\delta\subset T^\gamma$ and
$T^\delta_{[s]}\uparrow T^\gamma_{[s]}$ as $\delta\uparrow\gamma$. To make
this construction we merely have to choose a uniform variable $U(e)$ for
each edge $e$ in $T^\gamma$, and then for the construction of $T^\delta$
remove $e$ (and its successors) if and only if $U(e)>\delta/\gamma$. One
easily sees from the fact that the number of edges from $\la i_1,\dots,i_n\ra $
to $T^\gamma_{n+1}$ has a Poisson distribution with mean $\gamma$ that the
resulting tree has the same distribution as $T^\delta$.

Once we have (3.65) and (3.72) we obtain that
$$
\lim_{\delta\uparrow\gamma}\lim_{K\to\infty} R(\delta;\varepsilon,K)
=R^\varepsilon(T^\gamma)\quad\text{in distribution}.
$$ 
Finally, we may let $\varepsilon\downarrow 0$ by virtue of Proposition 1.
\qed\enddemo

\section {4. Some further remarks concerning the proof of Theorem 1} 
Theorem 1 in Grimmett and Kesten (1983) asserts that
$$
\gamma (n)R_n\to 2\bigg\{\int_{[0,\infty)}x^{-1}\,dF(x)\bigg\}^{-1}\quad
\text{ in probability }\tag 4.1 
$$
if $\gamma(n)\to\infty$ as $n\to\infty$. In view of Theorem 4 of Grimmett
and Kesten (1983) and a passage to subsequences we may restrict ourselves
to the case where
$$
\gamma (n)\to\infty\quad\text{but}\quad{\log\gamma(n)\over\log n}\to 0.\tag 4.2
$$
Also we only have to prove an upper bound for $R_n$, since Lemma 5 in
Grimmett and Kesten (1983) provides the necessary lower bound. 

To obtain
an upper bound for $R_n$, let $\eta>0$ be a small number and $K <\infty$
such that
$$
F(K)\ge 1-\eta.
$$
Define $m_n$ by (3.16) with $\gamma =\gamma (n)$ and set $s_n = (m_n)^{1/2}$.
One can then show that with probability tending to 1 (as 
$n\to\infty$)
the following events (4.3)--(4.5) do occur:
\roster
\item"{(4.3)}"  $\tau^0_{[m_n]}$ and $\tau^\infty_{[m_n]}$ as constructed in the
beginning of Sect.\ 3 contain subtrees $\widetilde T$ and $\wwtilde T$ with
roots at $0$ and $\infty$, respectively. $\widetilde T$ and $\wwtilde T$
are disjoint, and each has $m_n$ generations. Each vertex of $\widetilde
T$ which is not in the
$m_n$th generation of $\widetilde T$ has exactly
$\lfloor(1-2\eta)\gamma(n)\rfloor$ children in $\widetilde T$. The last
statement remains true when $\widetilde T$ is everywhere replaced by
$\wwtilde T$.

\item "{(4.4)}" Each edge in $\widetilde T$ not incident to $0$ and each
edge in $\wwtilde T$ not incident to $\infty$ has resistance at most $K$.

\item "{(4.5)}" For each vertex $x$ of the $s_n$th generation of
$\widetilde T$ and each vertex $y$ of the $s_n$th generation of $\wwtilde
T$ there exists a path from $x$ to $y$ of $(2m_n-2s_n+1)$ edges of
$K_{n+2}$ of resistance $\le K$. This path consists of three pieces: (i) a
path of $(m_n-s_n)$ edges in $\widetilde T$ in the tree of descendants of
$x$, going from $x$ to a vertex $u$ in the $m_n$th generation of
$\widetilde T$; (ii) a single edge from $u$ to a vertex $v$ of the
$m_n$th generation of $\wwtilde T$ (this edge does not belong to
$\widetilde T\cup\wwtilde T)$;\quad(iii) a path of $(m_n-s_n)$ edges in
$\wwtilde T$ from $v$ to $y$ in the tree of descendants of $y$. For
distinct pairs $x,y$ the edges between the $m_n$th generations of
$\widetilde{T}$ and $\wwtilde T$ in the corresponding paths are different.
\endroster

The proofs of (4.3)--(4.5) are analogous to those of Lemmas 7, 10 and 11.
Note that (4.3)--(4.5) contain no statement about the resistances of the
edges incident to $0$ in $\widetilde T$ and the edges
incident to $\infty$ in $\wwtilde T$, except that they are edges of
$\tau^0_{[m_n]}$ and $\tau^\infty_{[m_n]}$. Accordingly, the conditional
distribution of the resistances of these edges is simply $F$. Since we are
only looking for an upper bound for $R_n$ we may raise all other
resistances in $\widetilde T$ and $\wwtilde T$, and those between the
$m_n$th generation of $\widetilde T$ and the $m_n$th generation of
$\wwtilde T$ to $K$. From here on the proof proceeds exactly as in Lemmas
12 and 13 with the following replacements for $T'_{[s]}$, $T''_{[s]}$ and
$M(n)$. $T'_{[s]}$ and $T''_{[s]}$ each are trees --- rooted at $0$ and
$\infty$, respectively --- of $s=s_n$ generations, in which each vertex
except those of the $s_n$th generation has exactly $\lfloor
(1-2\eta)\gamma(n)\rfloor$ children. Each edge not incident to one of the
roots has resistance $K$, while all edges incident to one of the roots
have independent resistances, chosen according to the distribution
function $F$. $M(n)$ consists of $T'_{[s]}$,  $T''_{[s]}$ plus an edge
between each pair of vertices $x$, $y$ with $x$ in the $s$th generation of
$T'$ and $y$ in the $s$th generation of $T''$. The latter edges are
distinct and each has resistance
$$
2m_nK\lbrace\lfloor(1-2\eta)\gamma(n)\rfloor\rbrace^{s_n}. \tag 4.6
$$     
(4.6) takes the place of (3.31). With these replacements the proofs of
Lemma 12 and 13 need no significant changes. ${\Cal A}'_s$ (${\Cal
A}''_s$) is simply replaced by the full $s$th generation of
$T'$ ($T''$) and no appeal to Lemma 6 is needed, since now the analogue
of $\pi(\la x\ra ,T_{[m]},R,s)$ becomes
$$
\Pr\lbrace X.\text{ reaches } \la 0\ra  \text{ before } 
T_s\mid X_0=\la x\ra \rbrace \tag 4.7
$$
for $\la x\ra \in T_{s-1}$, and this last probability is always 
$\le (2L_n/(2L_n+K))^{s-1}$, where $2L_n$ is the resistance of an infinite 
tree in which each vertex has $\lfloor(1-2\eta)\gamma(n)\rfloor$ children 
and each edge has resistance $K$ (compare proof of Lemma 6; 
actually $L_n\to 0$ as $n\to\infty$, but that is not important). 
It is also easy to see that for each $K$
$$
\multline
(1-2\eta)\gamma (n)\{\text{resistance between } 0\text{ and }
 T'_s \text{ in }T'_{[s]}\}\to \left\lbrace\int {1\over x}\,dF(x)
\right\rbrace^{-1}\\
\text{in probability as } n\to\infty.
\endmultline
$$
It is not necessary to remove truncations --- as done in Lemma 14 --- in the
present case. In particular Prop.\ 1 is not needed when (4.2) holds. We
merely have to take the limit as $\eta\downarrow 0$ in the preceding
estimates.

\enddocument